\documentclass[%
11pt,
reprint,
onecolumn,
tightenlines,
superscriptaddress,
preprintnumbers,
nofootinbib,
amsmath,amssymb,amsthm,
physrev,
eqsecnum,tikz,
]{revtex4-2}

\usepackage{isomath}
\usepackage{amsmath,amsthm}
\usepackage{amsbsy}
\usepackage{amssymb}
\usepackage{amscd}
\usepackage{amsfonts}
\usepackage{stmaryrd}
\usepackage{siunitx}
\usepackage{euscript}
\usepackage[utf8]{inputenc}
\usepackage[T1]{fontenc}
\usepackage{newtxtext} 
\everymath{\displaystyle}
\usepackage{exscale}
\usepackage{microtype}
\usepackage{hyperref}
\usepackage{booktabs}
\usepackage{algorithm}
\usepackage{algpseudocode}
\usepackage{graphicx}
\usepackage{boxedminipage}
\usepackage{calc}
\usepackage[usenames,dvipsnames]{xcolor}
\graphicspath{ {media/} }
\usepackage[caption=false,justification=centerlast]{subfig}

\usepackage{setspace}
\usepackage{enumitem}
\setitemize{noitemsep,topsep=0pt,parsep=0pt,partopsep=0pt}
\setenumerate{noitemsep,topsep=0pt,parsep=0pt,partopsep=0pt}
\setdescription{noitemsep,topsep=0pt,parsep=0pt,partopsep=0pt}

\usepackage[colorinlistoftodos, color=green!40,prependcaption]{todonotes}
\setuptodonotes{inline}

\usepackage{soul} 
\usepackage[normalem]{ulem}

\usepackage{orcidlink}
\usepackage{siunitx}

\usepackage[small]{titlesec}

\titlespacing*{\section}{0pt}{12pt plus 4pt minus 2pt}{2pt plus 2pt minus 2pt}
\titlespacing*{\subsection}{0pt}{12pt plus 4pt minus 2pt}{2pt plus 2pt minus 2pt}
\titlespacing*\subsubsection{0pt}{12pt plus 4pt minus 2pt}{2pt plus 2pt minus 2pt}
\titlespacing*\paragraph{0pt}{12pt plus 4pt minus 2pt}{2pt plus 2pt minus 2pt}

\makeatletter
    \renewcommand*{\thesection}{\arabic{section}}
    \renewcommand*{\thesubsection}{\thesection.\Alph{subsection}}
    \renewcommand*{\p@subsection}{}
    \renewcommand*{\thesubsubsection}{\thesubsection.\arabic{subsubsection}}
    \renewcommand*{\p@subsubsection}{}
\makeatother

\setuptodonotes{inline}
\usepackage{isomath}
\usepackage{amsmath}
\usepackage{amssymb}
\usepackage{amscd}
\usepackage{amsfonts}

\DeclareMathOperator{\divergence}{div}
\DeclareMathOperator{\divS}{div_S}

\newcommand{\dm}{\ \mathrm{d}}

\usepackage[utf8]{inputenc}
\usepackage{amsmath}
\usepackage{amsthm}
\usepackage{amssymb}
\usepackage{bbm}
\usepackage{graphicx}
\usepackage{blindtext}
\usepackage{stmaryrd}
\usepackage{color}
\usepackage{lipsum}
\usepackage{tikz}
\usetikzlibrary{calc}
\usetikzlibrary{angles,quotes}
\usepackage[rightcaption]{sidecap}
\graphicspath{ {images/} }
\usepackage{dirtytalk}
\usepackage{accents}
\usepackage{contour} 
\usepackage{ulem}    
\usepackage{hyperref}
\usepackage{mathtools}

\newcommand{\Scr}[1]{\mathbb{#1}}
\newcommand{\scr}[1]{\mathcal{#1}}

\newcommand{\ve}[1]{\boldsymbol{#1}}

\newcommand{\utilde}[1]{\underaccent{\tilde}{#1}}

\contourlength{0.8pt}
\newcommand{\myuline}[1]{\uline{\phantom{#1}}\llap{\contour{white}{#1}}}

\newtheorem{theorem}{Theorem}[section]
\AtEndEnvironment{theorem}{\null\hfill\qedsymbol}

\newtheorem{lemma}[theorem]{Lemma}

\theoremstyle{remark}
\newtheorem{remark}{Remark}[section]
\AtEndEnvironment{remark}{\null\hfill\qedsymbol}

\theoremstyle{definition}

\AtEndEnvironment{definition}{\null\hfill\qedsymbol}

\newtheorem{mydef}{Definition}
\AtEndEnvironment{mydef}{\null\hfill\qedsymbol}

\newcommand{\defn}[2]{\begin{mydef}{\underline{#1}:}
#2
\end{mydef}}
\newcommand{\rmrk}[1]{\begin{remark}
#1
\end{remark}}

\theoremstyle{definition}

\AtEndEnvironment{exmp}{\null\hfill\qedsymbol}

\newcommand{\beol}[1]{\begin{equation}
#1
\end{equation}}

\newcommand{\bml}[1]{\begin{equation}
\begin{split}
#1
\end{split}
\end{equation}}

\graphicspath{{ImagesWithWords/}}

\begin{document}

\preprint{To appear in Mathematics and Mechanics of Solids (DOI: \href{https://doi.org/10.1177/10812865251361554}{10.1177/10812865251361554})}

\title{Nonlocal Electrostatics and Boundary Charges in Continuum Limits of Two-Dimensional Materials}

\author{Shoham Sen}
    \email{shohams@andrew.cmu.edu}
    \affiliation{Department of Civil and Environmental Engineering, Carnegie Mellon University}

\author{Yang Wang}%
    \affiliation{Pittsburgh Supercomputing Center}%

\author{Timothy Breitzman}
    \affiliation{Air Force Research Laboratory}
    
\author{Kaushik Dayal}
    \affiliation{Department of Civil and Environmental Engineering, Carnegie Mellon University}
    \affiliation{Center for Nonlinear Analysis, Department of Mathematical Sciences, Carnegie Mellon University}
    \affiliation{Department of Mechanical Engineering, Carnegie Mellon University}
    
\date{\today}

\begin{abstract}
    Two-dimensional (2D) electronic materials are of significant technological interest due to their exceptional properties and broad applicability in engineering. 
    The transition from nanoscale physics, that dictates their stable configurations, to macroscopic engineering applications requires the use of multiscale methods to systematically capture their electronic properties at larger scales.
    A key challenge in coarse-graining is the rapid and near-periodic variation of the charge density, which exhibits significant spatial oscillations at the atomic scale.
    Therefore, the polarization density field --- the first moment of the charge density over the periodic unit cell --- is used as a multiscale mediator that enables efficient coarse-graining by exploiting the almost-periodic nature of the variation.
    Unlike the highly oscillatory charge density, the polarization varies over lengthscales that are much larger than the atomic,  making it suitable for continuum modeling.
    
    In this paper, we investigate the electrostatic potential arising from the charge distribution of 2D materials. 
    Specifically, we consider a sequence of problems wherein the underlying lattice spacing vanishes and thus obtain the continuum limit.
    We consider 3 distinct limits: where the thickness is much smaller than, comparable to, and much larger than the in-plane lattice spacing.
    These limiting procedures provide the homogenized potential expressed in terms of the boundary charge and dipole distribution, subject to the appropriate boundary conditions that are also obtained through the limit process.
    Further, we demonstrate that despite the intrinsic non-uniqueness in the definition of polarization, accounting for the boundary charges ensures that the total electrostatic potential, the associated electric field, and the corresponding energy of the homogenized system are uniquely determined.
\end{abstract}

\maketitle



\section{Introduction}

Two-dimensional (2D) materials have been the focus of significant interest since the isolation of graphene in 2004 \cite{gerstner2010nobel}. 
These materials exhibit unique electronic, mechanical, and optical properties, making them highly promising for a wide range of technological applications. 
In particular, 2D materials are being explored for energy-harvesting devices such as flexoelectric and piezoelectric nanogenerators \cite{wu2014piezoelectricity,lee2017reliable}, triboelectric nanogenerators \cite{zi2015triboelectric,pace20232d,ghosh2022ferroelectricity}, actuation systems \cite{wang2016subatomic}, piezotronics \cite{wu2013taxel,xue2015influence,pradel2013piezotronic,zhang2021piezotronics}, and non-volatile random-access memory \cite{wang20232d}. 
Moreover, these materials hold great promise for exhibiting strong flexoelectric responses due to their ability to sustain large strain gradients through bending deformations \cite{liang2024advancements,surmenev2023influence}. 
These materials have been utilized in a diverse range of applications, including but not limited to flexible electronics, optoelectronic devices, field-effect transistors (FETs), sensors, ultrafast lasers, batteries, and supercapacitors \cite{glavin2020emerging, kumbhakar2023prospective,liu2016localized,lin2023recent}. 
A  notable aspect of two-dimensional (2D) materials is their ability to exhibit large strain gradients, which can induce significant flexoelectric responses. From an application perspective, these materials hold considerable promise for the development of advanced sensors and energy harvesting devices.

The functional response of these materials (also for bulk 3D materials) at the continuum lengthscale is typically described by the polarization, which serves as a multiscale mediator that captures key aspects of the atomic-scale charge distribution
\cite{james1994internal,xiao2005influence,jha2023discrete,jha2023atomic,tadmor-EffHamil,grasinger2021flexoelectricity,grasinger2021architected,grasinger2020statistical,marshall2014atomistic,cicalese2009discrete,bach2020discrete,alicandro2008continuum,sharp1994electrostatic,grasinger2020statistical, grasinger2021architected,grasinger2022statistical,sen2024nonuniqueness}. 
The role of polarization in piezoelectric and ferroelectric phenomena has been the focus of extensive research: e.g., 
in bulk ferroelectrics \cite{indergand2021effect,cheng2024growth,indergand2020phase,kannan2022kinetics,wojnar2017linking}; 
in fundamental studies of polarization \cite{tagantsev1987pyroelectric,tagantsev1986piezoelectricity,tagantsev2013origin};
in phase-transforming settings \cite{bucsek2019direct,jalan2010large,nunn2018frequency,sen2024multifunctional,james2021conversion,chen2011weak};
in functional composites \cite{nepal2020toward,biswas2013coherent,nepal2012high};
through continuum phase-field modeling \cite{ross2023thermodynamics,wang2021inverse,zorn2022q,chen2008phase,li2016origin,zhoudomain,pohlmann2017thermodynamic,wang2024theory};
in low-dimensional nanostructures \cite{zhang2019room,tian2020hexagonal,jha2023atomic,jha2023discrete};
in biomembranes \cite{steigmann2018mechanics,mathew2024electro,mathew2024active};
and flexoelectricity \cite{liu2013flexoelectricity, wang2019flexoelectricity, krichen2016flexoelectricity,abdollahi2014computational, Irene1, codony2021transversal}.

An important issue is that the polarization, defined as the dipole moment per unit cell, is inherently non-unique due to the arbitrary choice of the unit cell. 
Consequently, different selections of the unit cell yield distinct dipole moments (including reversing direction). 
Since polarization serves as a fundamental multiscale descriptor of the electronic behavior in continuum models, its non-uniqueness poses challenges to the consistency of continuum models. Various approaches have been proposed to address this issue. A widely adopted framework within the physics community is the \textit{Modern Theory of Polarization} \cite{Resta-Vanderbilt,vanderbilt_2018,king1993theory,Vanderbilt_surface}, which introduces an alternative definition of polarization through the change in the Berry phase of the wavefunction during an adiabatic process. The central assertion is that, while polarization itself is non-unique, the change in polarization is uniquely defined modulo a polarization quantum. The assumption of adiabaticity implies that eigenstates remain unchanged throughout the transition, allowing one to obtain an expression for the electric current using the eigenstates at equilibrium. This approach circumvents the issue of non-uniqueness inherent in defining polarization as the dipole moment per unit cell. By expressing polarization in terms of its change rather than its absolute value, this approach ensures well-defined values modulo the polarization quantum. An alternative approach, popular in the continuum mechanics community, defines polarization as the energy conjugate to the electric field \cite{toupin1956elastic}. Both, by definition, yield unique values.

In prior work \cite{Sen2022,sen2024rigorous,sen2024nonuniqueness}, we addressed the non-uniqueness in polarization from a classical perspective (dipole moment per unit cell definition). We demonstrated that for every choice of polarization, there exists a corresponding surface charge density. Different choices of polarization in the bulk produce different choices of surface charge density on the boundary. Together, these quantities compensate each other to yield unique potentials, fields, and energies. We showed in \cite{sen2024nonuniqueness} how the different definitions of polarization were connected. 

\paragraph*{Contributions of this work.}

This work aims to derive a series of polarization theorems for two-dimensional (2D) materials. 
We consider the electrostatic (Poisson) equation for a charge distribution supported on a thin film, and analyze the limit of vanishing thickness and lattice spacing. 
Previously, we addressed the non-uniqueness of polarization for 3D materials \cite{sen2024nonuniqueness,Sen2022}, by examining a sequence of Poisson equations associated with a corresponding sequence of charge distributions defined on a lattice. 
By appropriately scaling the charge distribution, we derived a homogenized Poisson equation posed with respect to the polarization. 
Analyzing this limit provides a resolution to the issue of polarization non-uniqueness.

In this paper, we follow a similar approach to investigate the non-uniqueness of polarization. Given two scaling parameters --- $l$ for the lattice and $h$ for the thickness --- we explore different asymptotic regimes:  
(i) $l \ll h \to 0$,  
(ii) $l \approx h \to 0$, and  
(iii) $h \ll l \to 0$.  
We demonstrate that, in the limit, the surface charge compensates for the tangential polarization \footnote{In the more general case where free charge is present (lack of charge-neutral unit cells), the free charge serves to compensate for the non-uniqueness in polarization.}, while the normal polarization gives rise to a double-layer potential.
The homogenized or coarse-grained problem enables us to replace the full problem real problem with its homogenized limit which is significantly more tractable for analysis and numerical approximation given that the polarization varies on scales much larger than atomic.

The sequence of Poisson equations that we analyze are:
\bml{-\Delta\Phi_{l,h}&=\rho_{l,h}\mathbbm{1}_{\Omega\times(-h,h)}~~\text{in}~\Scr{R}^3\\
\Phi_{l,h}(\ve{r})&=0~~\text{as}~|\ve{r}|\to\infty,\label{Problem-Intro}}
where the charge distribution is parameterized by $l$ and $h$. We are interested in the limit that $l,h\to0$. Since these parameters can vary independently, we aim to investigate different asymptotic regimes based on the relative scaling of \( l \) and \( h \).

In the limit \( l \ll h \to 0 \), with \( \rho_{l,h} = \frac{\tilde{\rho}_{l,h}}{lh} \), the asymptotic behavior of \eqref{Problem-Intro} is characterized by the following limit:
       \bml{-\Delta\Phi_{0}&=\left({q}-\dfrac{\divS\left(J_0 {\ve{p}}_p\right)}{J_0}\right)\mathbbm{1}_{\Omega}+\mathbbm{1}_{\partial\Omega}\left( {\sigma}+ {\ve{p}}_p\cdot \ve{n}\right)~~\text{in}~\Scr{R}^3\\
        \Phi_{0}(\ve{r})&=0~~\text{as}~|\ve{r}|\to\infty.\label{PDE-Case1-Intro}}
        where $\divS$ is the surface divergence, $J_0$ is the surface Jacobian associated with the manifold $\Omega$, $\ve{p}_p$ is the tangential or planar polarization, $q$ is the free charge, $\sigma$ is the surface charge distribution and $\ve{n}$ is normal to $\partial\Omega$ and tangent to $\Omega$.

In the limit \( h \ll l \to 0 \), with \( \rho_{l,h} = \frac{\tilde{\rho}_{l,h}}{h^2} \), the asymptotic behavior of \eqref{Problem-Intro} is characterized by the following limit:
\bml{-\Delta\Phi_{0}&= q+{p}_3\dfrac{\partial \mathbbm{1}_{\Omega}}{\partial \nu}~~\text{in}~\Scr{R}^3\\
        \Phi_{0}(\ve{r})&=0~~\text{as}~|\ve{r}|\to\infty.\label{PDE-Case3-Intro}}
        where $q$ is the free charge and ${p}_3$ is the normal polarization. Here $\ve{\nu}$ is the direction normal to the manifold. 

Finally, in the limit \( h \approx \alpha l \to 0 \), with \( \rho_l = \frac{\tilde{\rho}_l}{l^2} \), the corresponding limit of \eqref{Problem-Intro} is given by:
\bml{-\Delta\Phi_{0}&=\alpha\left({q}-\dfrac{\divS\left(J_0 {\ve{p}}_p\right)}{J_0}\right)\mathbbm{1}_{\Omega}+\alpha^2p_3\dfrac{\partial\mathbbm{1}_{\Omega}}{\partial\nu}+\mathbbm{1}_{\partial\Omega}\alpha\left( {\sigma}+ {\ve{p}}_p\cdot \ve{n}\right)~~\text{in}~\Scr{R}^3\\
        \Phi_{0}(\ve{r})&=0~~\text{as}~|\ve{r}|\to\infty.\label{PDE-Case2-Intro}}
        where $\divS$ is the surface divergence, $J_0$ is the surface Jacobian associated with the manifold $\Omega$, $\ve{p}_p$ is the tangential polarization, ${p}_3$ is the normal polarization, $q$ is the free charge, $\sigma$ is the surface charge, $\ve{n}$ is normal to $\partial\Omega$ and tangent to $\Omega$ and $\ve{\nu}$ is the normal to $\Omega$.

The results can be understood by observing that the normal component of the polarization remains uniquely defined and contributes to the formation of a double-layer potential, whereas the tangential component, along with the free charge, compensates for the surface charge distribution. 
Free charge typically arises in scenarios where there is a deviation from charge-neutral unit cells. 
Consequently, in cases where no free charge is present, the surface charge compensates for the non-uniqueness of polarization in the bulk. 
This distinction arises from the inherent non-uniqueness in the choice of the unit cell within the plane.

\paragraph*{Organization.}

This paper is organized as follows. In Section \ref{Notation}, we summarize the notation and fundamental definitions used throughout the paper. In Section \ref{Heuristic}, we present a heuristic understanding of the results in this paper. The problem is formulated in Section \ref{PF}. Section \ref{Homo} presents the homogenization theorem \ref{Theorem-Homo} that serves as a key component in the derivation of our main results. The proof of theorem \ref{Theorem-Homo} is divided into three parts one for each asymptotic regime. We address $l\ll h$ in Section \ref{Proof1}, $l\approx h$ in \ref{Proof2} and $h\ll l$ in Section \ref{Proof3}. We discuss the resolution of the non-uniqueness of the polarization field in Section \ref{Non-Uniqueness}. 
Appendix \ref{Appendix} provides the derivation of the double layer potential.

\section{Notation}\label{Notation}
We first present a summary of all the notations used throughout the paper for the readers to refer to:
\begin{itemize}
\item{\makebox[2cm]{$\Scr{Z}$\hfill} The set of integers.}
\item{\makebox[2cm]{$\Scr{R}$\hfill} The set of real numbers. Further, $\Scr{R}_{>0}$ and $\Scr{R}_{<0}$ represent the strictly positive and strictly negative subset.
}
\item{\makebox[2cm]{$\Scr{R}^n$\hfill} $\Scr{R}\times\Scr{R}\times \Scr{R}\times \ldots\times \Scr{R}$. The cartesian product of $n$ copies of the set of real numbers $\Scr{R}$.}
\item{\makebox[2cm]{$\scr{C}^r(X,Y)$\hfill} The space of continuous functions $f:X\to Y$ with continuous first $r$ derivatives.}
\item{\makebox[2cm]{$\mathcal{C}_0(X,Y)$\hfill} The space of compactly supported continuous functions $f:X\to Y$.}
\item{\makebox[2cm]{$\{\}$\hfill} The null set.}
\item{\makebox[2cm]{$P\oplus Q$\hfill} Minkowski sum of two sets $P,Q\subset\Scr{R}^n$ as $\{w\in\Scr{R}^n|\exists p\in P,q\in Q \text{ such that } w=p+q\}$.}
\item{\makebox[2cm]{$B^n_{r}( {x})$\hfill} Ball of radius $r$ centered at $ {x}$ in dimension $n$.}
\item{\makebox[2cm]{$l\Box$\hfill} Denotes scaled unit cell of dimensions $l\mathsf{x}l$.}
\item{\makebox[2cm]{$\Box\equiv l\Box\big\vert_{l=1}$\hfill} A unit cell.}
\item{\makebox[2cm]{$T(l\Box)$\hfill} The set of all $l\times l$ cells contained in $T\subset\Scr{R}^2$.}
\item{\makebox[2cm]{$l\Delta$\hfill} Denotes scaled partial unit cell.}
\item{\makebox[2cm]{$\Delta\equiv l\Delta\big\vert_{l=1}$\hfill} A rescaled partial unit cell.}
\item{\makebox[2cm]{$T(l\Delta)$\hfill} The set of all $l\Delta$ in $T\subset\Scr{R}^2$.}
\item{\makebox[2cm]{$\psi$\hfill} Denotes the parametric map.}
\item{\makebox[2cm]{$\scr{L}$\hfill} A lattice constructed from the basis lattice vectors $\{\ve{e}_i\}_{i=1}^{2}$.}
\item{\makebox[2cm]{$l\scr{L}$\hfill} Lattice constructed from lattice vectors $\{l\ve{e}_i\}_1^2$.}
\item{\makebox[2cm]{$\scr{I}$\hfill} An indexing set for the scaling parameter $l$ (or $h$), with $\inf \scr{I} = 0$.}
\item{\makebox[2cm]{$|\cdot|$\hfill} The notation is context-dependent and represents distinct mathematical concepts accordingly (Cardinality of a set, Lebesgue measure of a set, and absolute value).}
\end{itemize}

In this paper, vectors in \(\Scr{R}^2\) and \(\Scr{R}^3\) will be represented using boldface notation. We frequently consider the decomposition \(\Scr{R}^3 = \Scr{R}^2 \times \Scr{R}\) to distinguish between in-plane and out-of-plane components. The in-plane (tangential) components will be indicated by the subscript \( p \), while the out-of-plane component will be denoted by the subscript \( 3 \). To maintain clarity, we refrain from using boldface for vectors in general spaces. 

To rigorously define the two-dimensional (2D) manifold \(\Omega\) and the charge distribution localized in its vicinity, we consider \(\Omega\) as a submanifold of \(\mathbb{R}^3\). There exists a parametric mapping \(\psi: \mathbb{R}^3 \to \mathbb{R}^3\) satisfying \(\psi(T \times (-h,h)) = \Omega \times (-h,h)\), where \(T \subset \mathbb{R}^2\) represents a reference domain in the parameter space. We denote vectors in the physical space \(\Scr{R}^3\) by \(\ve{r}\) and those in the parameterization space \(\Scr{R}^3\) by \(\ve{x}\), such that the relation \(\ve{r}=\psi(\ve{x})\) holds. Additionally, we impose \(\psi(T \times \{0\}) = \Omega\), ensuring that the reference configuration is mapped onto the manifold. To formalize this further, we define \(\psi_0 = \psi \vert_{x_3=0}\) as the restriction of \(\psi\) to the plane \(x_3=0\). The Jacobian associated with $\psi$ will be denoted as $J$, while the Jacobian associated with $\psi_0$ will be denoted $J_0$. We denote $\ve{n}$ as the tangent vector on $\Omega$, normal to $\partial\Omega$, while $\ve{\nu}$ is the normal to $\Omega$.

The surface divergence on the manifold \(\Omega\) is denoted by \(\divS\), while the divergence in the parameter space is denoted by \(\divergence\). Additionally, the notation \(\divergence_p\) is employed to indicate the divergence restricted to the reference domain \(T\). Subscripts are explicitly included when differentiating the Green's function to clarify the variable with respect to which differentiation is performed, given its dependence on two arguments. In particular, the subscript \(\ve{x}'\) designates the three-dimensional gradient with respect to the second argument of the Green's function, whereas the subscript \(\ve{x}'_p\) denotes the projection of this gradient onto the plane \(x_3 = 0\).

To facilitate the derivation of the polarization theorem, we introduce some notations. Prior to this, we provide a concise review of the relevant mathematical definitions necessary for the subsequent analysis.

\begin{defn}
{\textbf{$\scr{C}^r$-Diffeomorphism}} {Let $U, V \subset \mathbb{R}^n$ be open sets. A map $\phi \in \scr{C}^r(U, V)$ is called a $\scr{C}^r$-diffeomorphism iff:
\begin{enumerate}
    \item[i)] $\phi$ is bijective,
    \item[ii)] $\phi^{-1} \in \scr{C}^r(V, U)$.
\end{enumerate}}
\end{defn}

\begin{defn}
{\textbf{$k$-dimensional submanifold of $\mathbb{R}^n$}} {Let $1 \leq k \leq n-1$. A subset $M \subset \mathbb{R}^n$ is a $k$-dimensional submanifold of $\mathbb{R}^n$ if, for each $x \in M$, there exists a neighborhood $U$ of $x$ in $\mathbb{R}^n$ and a map $\phi: U \to \mathbb{R}^n$ such that:
\begin{enumerate}
    \item[i)] $\phi: U \to \phi(U)$ is a diffeomorphism,
    \item[ii)] $\phi(M \cap U) = \phi(U) \cap \{x \in \mathbb{R}^n \mid x_{k+1} = \cdots = x_n = 0\}$.
\end{enumerate}}
\end{defn}

\rmrk{Intuitively speaking, a manifold is a set that can be locally flattened. Practitioners of continuum mechanics may think of a diffeomorphism as the deformation map that maps the reference configuration to the deformed configuration. In this paper, we are interested in 2D manifolds, so one may think of a 2D manifold as the deformed configuration of a flat reference surface.}

\begin{defn}
{\textbf{Local Parameterization}} {Let $M$ be a submanifold of $\mathbb{R}^n$ and $y \in M$. A local parameterization of $M$ at $y$ is a map $\psi: V \to \mathbb{R}^n$, with $\psi(x) \in M$ for all $x \in V$, where $V \subset \mathbb{R}^k$ is open. The parameterization must satisfy:
\begin{enumerate}
    \item[i)] $0 \in V$ and $\psi(0) = y$,
    \item[ii)] There exists a neighborhood $U$ of $y$ in $\mathbb{R}^n$ such that $\psi: V \to U \cap M$ is bijective,
    \item[iii)] $\psi^{-1}: U \cap M \to V$ is continuous,
    \item[iv)] The derivative $D\psi(x)$ has full rank $k$ for all $x \in V$.
\end{enumerate}}
\end{defn}

\begin{theorem}\label{MST} \myuline{\textbf{Multivariable substitution theorem}}\\
Let $M$ be a k-dimensional submanifold on $\Scr{R}^n$. $V\subset\Scr{R}^k$ and $U\subset\Scr{R}^n$ open with $\psi:V\to M\cap U$ a parameterization of $M$. If $f:M\to\Scr{R}$ is a continuous function such that $\text{supp}(f)\subset U$, then we have
\beol{\int_M f \dm S=\int_V(f\circ\phi)(x)\sqrt{\text{det}\left(D_x\psi^T D_x\psi\right)}\dm x,}
where \(\dm S\) denotes the surface measure, following the convention that integrals over \(n\)-dimensional spaces are typically regarded as volume integrals, while those over lower-dimensional subsets are interpreted as surface integrals.
\end{theorem}

We now provide a brief introduction to the characteristic function. These will be used in writing the homogenized solutions. Given any domain $\Omega \subset \Scr{R}^3$, the characteristic function $\mathbbm{1}_{\Omega}(\ve{x})$ is defined as:
\begin{equation}
    \mathbbm{1}_{\Omega}(\ve{x}) = 
    \begin{cases}
        1, & \ve{x} \in \Omega, \\
        0, & \ve{x} \notin \Omega.
    \end{cases}
\end{equation}

We now introduce the distributional derivative of the characteristic function:
\begin{equation}
    \int_{\Scr{R}^3} f(\ve{x}) \dfrac{\partial \mathbbm{1}_{\Omega}(\ve{x})}{\partial \nu} \, \mathrm{d} \ve{x} 
    = -\int_{\Scr{R}^3} \mathbbm{1}_{\Omega}(\ve{x}) \dfrac{\partial f(\ve{x})}{\partial \nu} \, \mathrm{d} \ve{x} 
    = -\int_{\Omega} \dfrac{\partial f(\ve{x})}{\partial \nu} \, \mathrm{d} \ve{x}, 
    \quad \forall f \in \mathcal{C}^{\infty}_0(\Scr{R}^3),
\end{equation}
where \(\ve{\nu}:\Omega\to\Scr{R}^3\) denotes a vector field. 

For the purpose of our analysis, we will set \(\ve{\nu}\) as the unit normal vector to a manifold \(\Omega\subset\Scr{R}^3\). To this end, we define,
\beol{\dfrac{\partial \mathbbm{1}_{\Omega}(\ve{x})}{\partial \nu}=\lim_{t\to0}\left(\dfrac{\mathbbm{1}_{\Omega}(\ve{x})-\mathbbm{1}_{\Omega_t}(\ve{x})}{t}\right),}
where $\Omega_t=\{\ve{x}-\ve{\nu}t\vert\ve{x}\in\Omega\}$.

To rigorously define a polarization distribution, it is essential to identify regions over which homogenization can be performed. Since the charge distribution undergoes scaling, these regions must be scaled accordingly. A natural and mathematically consistent choice is to consider regions forming a periodic tessellation of \( T\subset\mathbb{R}^2 \), as this ensures that the scaling behavior follows directly from their geometric structure. The construction of a periodic tessellation requires three fundamental components: a reference origin, a periodic lattice, and a fundamental domain.

Having chosen an origin for our coordinate system, we define a lattice $\scr{L}$ using two arbitrary basis vectors $\ve{e}_1, \ve{e}_2 \in \mathbb{R}^2$, which generate our lattice on $T$:
\beol{\scr{L} := \left\{ \sum_1^2 n^i \ve{e}_i \mid n^1, n^2 \in \mathbb{Z} \right\}.}

A lattice structure \( \psi(\scr{L}) \) on the manifold \( \Omega \) can be constructed by applying the mapping \( \psi \) to the lattice \( \scr{L} \). The scaled version of this lattice, obtained by introducing a scaling factor \( l \in (0,1] \), is given by:
\beol{l\scr{L} := \left\{ \sum_1^2 n^i l \ve{e}_i \mid n^1, n^2 \in \mathbb{Z} \right\}.}
The corresponding scaled lattice $\psi(l \scr{L})$ is obtained accordingly.

Subsequently, we introduce the set of admissible unit cells and their associated corner mappings, collectively denoted as \( \scr{AC} \). To this end, we select a reference unit cell \( {\Box}^{\prime} \subset \mathbb{R}^2 \) and a vector \( \ve{f} \in \mathbb{R}^2 \) satisfying the following conditions:
\begin{enumerate}
    \item ${\Box}^{\prime}$ is measurable in the sense of 2D Lebesgue measures,
    \item $\mathbb{R}^2 = \bigcup_{\ve{i} \in l\scr{L}} \ve{i} \oplus (l\ve{f} \oplus l{\Box}^{\prime})$,
    \item For all $\ve{i} \neq \ve{j} \in l\scr{L}$, $\ve{i} \oplus (l\ve{f} \oplus l{\Box}^{\prime}) \cap \ve{j} \oplus (l\ve{f} \oplus l{\Box}^{\prime}) = \{\}$,
\end{enumerate}
where the first condition guarantees the well-posedness of the polarization and charge distribution, while the subsequent two conditions ensure that the unit cell forms a periodic tessellation of $\Scr{R}^2$. In the above, $\oplus$ denotes the Minkowski sum (see notation summary list at the beginning of this section).

For the unit cell and corner map, we choose \((\Box^{\prime}, \ve{f}) \in \scr{AC}\) and define the scaled unit cell as \(l\Box := l\ve{f} \oplus l{\Box}^{\prime}\), with the corresponding corner map given by \(\hat{\ve{x}}_{p,l} := \ve{i} + l\ve{f}\), where \(\ve{i} \in l\scr{L}\). The subscript $p$ is used to denote planar, ie, $\ve{x}_p\in T$. The corner map establishes an association between each position vector and a reference point within the unit cell to which it belongs. For clarity, the subscript \(l\) will be omitted when there is no risk of ambiguity. The corresponding unit cells on the manifold are defined as $\psi(l\Box)$. Note that $\psi(l\Box)$ is dependent on the mapping $\psi$, as the positioning of the unit cell dictates its image under $\psi$.

Partial unit cells $l\Delta$ on $T$ are unit cells which are not entirely contained in $T$:
\beol{\hat{\ve{x}}_p \oplus l\Delta := \left\{\hat{\ve{x}}_p \oplus l\Box \cap T \mid \hat{\ve{x}}_p \in l\scr{L} \text{ and } \hat{\ve{x}}_p \oplus l\Box \cap T^c \neq \{\} \right\}.}

We define $T(l\Box)$ as the collection of full unit cells:
\[
T(l\Box) := \left\{ \hat{\ve{x}}_p \oplus l\Box \mid \hat{\ve{x}}_p \in l\scr{L}, \hat{\ve{x}}_p \oplus l\Box \subset T \right\}.
\]
Similarly, the collection of partial unit cells is denoted by $T(l\Delta)$ and defined as:
\[
T(l\Delta) := \left\{ \hat{\ve{x}}_p \oplus l\Box \cap T \mid \hat{\ve{x}}_p \oplus l\Box \cap T^c \neq \{\} \right\}.
\]

We can define $\Omega(l\Box) = \psi(T(l\Box))$ and $\Omega(l\Delta) = \psi(T(l\Delta))$.

\begin{rmrk}
{Each of the terms defined above holds independently of the parameterization chosen for $\Omega$. If two parametric maps $\psi: T \to \Omega$ and $\psi': T' \to \Omega$ exist, there also exist lattices $l\scr{L}$ and $l\scr{L}'$ such that $\psi(T \cap l\scr{L}) = \psi'(T' \cap l\scr{L}')$.}
\end{rmrk}

\section{Heuristics}\label{Heuristic}
In this section, we analyze the case of a flat manifold, which offers the advantage of eliminating curvature effects, thereby enabling a more tractable and straightforward analysis. Moreover, we employ a less rigorous approach by deriving several key results, including scaling behavior, through the application of a Taylor series expansion. Given that we are working with a flat manifold, we do not differentiate between the parameterization space and the real space. Additionally, due to the flat nature of the manifold, unit cells within the bulk are inherently charge neutral. Consequently, we do not consider the presence of free charge in this section.

Consider a sequence of charge distributions \(\{\rho_{l,h}\}_{l,h}:\Omega_{h}\to\Scr{R}\), where \(l\) and \(h\) are parameters that go to zero. Here, \(\Omega_h = \Omega \times (-h, h)\) denotes the support of the charge distribution, with \(\Omega \subset \Scr{R}^2\). The charge distribution is assumed to be locally \(l\)-periodic, with the additional condition that it remains charge-neutral within each unit cell. Our objective is to investigate the thin-film limit of the potential generated by these charge distributions. Given that the problem involves two scaling parameters, we will analyze different asymptotic regimes for the charge density, focusing on the relationship between \(l\) and \(h\). 

The potential generated by this charge distribution is given by 
\bml{-\Delta\Phi_{l,h}&=\rho_{l,h}~~\text{in}~\Scr{R}^3\\
\Phi_{l,h}&=0~~\text{as}~|x|\to\infty,\label{Problem}}

The solution to \eqref{Problem} is given by the convolution against the Green's function, 
\beol{\Phi_{l,h}(\ve{x})=\int_{\Omega_h}G(\ve{x},\ve{x}')\rho_{l,h}(\ve{x}')\dm\ve{x}'\label{solution}}

We proceed by decomposing the integral in \eqref{solution} into a sum of integrals over individual unit cells, utilizing the corner map to facilitate this process. Specifically, we employ the substitution \(\ve{x}' = \hat{\ve{x}}'_p + l\ve{y}' + hz'\ve{e}_3\), where \(\hat{\ve{x}}'_p \in \Omega\), \(\ve{y}' \in \Box\), and \(\ve{e}_3\) denotes the unit vector perpendicular to both \(\hat{\ve{x}}'_p\) and \(\ve{y}'\). 
\bml{\Phi_{l,h}(\ve{x})&=\sum_{\hat{\ve{x}}'_p\in\Omega(\Box)}l^2h\int_{\Box}\int_{-1}^1G(\ve{x},\hat{\ve{x}}'_p+l\ve{y}'+hz'\ve{e}_3)\rho_{l,h}(\hat{\ve{x}}'_p+l\ve{y}'+hz'\ve{e}_3)\dm \ve{y}'\dm z'\\
&+\sum_{\hat{\ve{x}}'_p\in\Omega(\Delta)}l^2h\int_{\Delta}\int_{-1}^1G(\ve{x},\hat{\ve{x}}'_p+l\ve{y}'+hz'\ve{e}_3)\rho_{l,h}(\hat{\ve{x}}'_p+l\ve{y}'+hz'\ve{e}_3)\dm \ve{y}'\dm z'\label{Sol-Expand}}

We are interested in the potential away from the domain. Thus, the Green's function should be smooth there. Taylor expanding the Green's function
\beol{G(\ve{x},\hat{\ve{x}}'_p+l\ve{y}'+hz'\ve{e}_3)=G(\ve{x},\hat{\ve{x}}'_p)+\left(l\ve{y}'+hz'\ve{e}_3\right)\cdot \nabla_{\ve{x}'} G(\ve{x},\hat{\ve{x}}'_p)+\scr{O}(l^2,lh,h^2).\label{Taylor-Series}}

Substituting \eqref{Taylor-Series} into \eqref{Sol-Expand},
\bml{\Phi_{l,h}(\ve{x})&=\sum_{\hat{\ve{x}}'_p\in\Omega(\Box)}l^2h G(\ve{x},\hat{\ve{x}}'_p)\int_{\Box}\int_{-1}^1\rho_{l,h}(\hat{\ve{x}}'_p+l\ve{y}'+hz'\ve{e}_3)\dm \ve{y}'\dm z'\\
&+\sum_{\hat{\ve{x}}'_p\in\Omega(\Box)}l^2h \nabla_{\ve{x}'} G(\ve{x},\hat{\ve{x}}'_p)\cdot\int_{\Box}\int_{-1}^1\left(l\ve{y}'+hz'\ve{e}_3\right)\rho_{l,h}(\hat{\ve{x}}'_p+l\ve{y}'+hz'\ve{e}_3)\dm \ve{y}'\dm z'\\
&+\sum_{\hat{\ve{x}}'_p\in\Omega(\Delta)}l^2h G(\ve{x},\hat{\ve{x}}'_p)\int_{\Delta}\int_{-1}^1\rho_{l,h}(\hat{\ve{x}}'_p+l\ve{y}'+hz'\ve{e}_3)\dm \ve{y}'\dm z'}

Now using charge neutrality, we can get rid of the first term. We define the inplane and out of plane dipole moment as
\bml{\Phi_{l,h}(\ve{x})&=\sum_{\hat{\ve{x}}'_p\in\Omega(\Box)}l^2 \nabla_{\ve{x}'} G(\ve{x},\hat{\ve{x}}'_p)\cdot\left(hl\ve{p}_{p,l,h}(\hat{\ve{x}}'_p)+h^2p_{3,l,h}(\hat{\ve{x}}'_p)\ve{e}_3\right)+\sum_{\hat{\ve{x}}'_p\in\Omega(\Delta)}l G(\ve{x},\hat{\ve{x}}'_p)lh\sigma_{l,h}(\hat{\ve{x}}'_p),}
where we defined
\bml{\ve{p}_{p,l,h}(\hat{\ve{x}}_p')&=\int_{\Box}\int_{-1}^1\ve{y}\rho_{l,h}(\hat{\ve{x}}'_p+l\ve{y}'+hz'\ve{e}_3)\dm \ve{y}'\dm z'\\
{p}_{3,l,h}(\hat{\ve{x}}_p')&=\int_{\Box}\int_{-1}^1z\rho_{l,h}(\hat{\ve{x}}'_p+l\ve{y}'+hz'\ve{e}_3)\dm \ve{y}'\dm z'\\
\sigma_{l,h}(\hat{\ve{x}}_p')&=\int_{\Delta}\int_{-1}^1\rho_{l,h}(\hat{\ve{x}}'_p+l\ve{y}'+hz'\ve{e}_3)\dm \ve{y}'\dm z'\label{pol-charge}}

We can approximate $\textstyle\sum_{\hat{\ve{x}}'_p\in\Omega(\Box)}l^2$ with $\textstyle\int_{\Omega}\dm \ve{x}_p'$ and $\textstyle\sum_{\hat{\ve{x}}'_p\in\Omega(\Delta)}l$ with $\textstyle\int_{\partial\Omega}\dm S_{\ve{x}_p'}$ where $\hat{\ve{x}}_p\mapsto\ve{x}_p$ in the limit of $l\to0$. We can thus approximate the answer with this,
\bml{\Phi_{l,h}(\ve{x})&\approx\int_{\Omega}\dm\ve{x}_p' \nabla_{\ve{x}'} G(\ve{x}, {\ve{x}}'_p)\cdot\left(hl\ve{p}_{p,l,h}(\ve{x}'_p)+h^2p_{3,l,h}(\ve{x}'_p)\ve{e}_3\right)+\int_{\partial\Omega}\dm S_{\ve{x}_p'} G(\ve{x}, {\ve{x}}'_p)lh\sigma_{l,h}(\ve{x}'_p),\label{approx-sol}}

\subsection{Case I. $h\ll l\to0$}\label{Case I}
For this limit to make sense, we consider the scaling, $\rho_{l,h}=\dfrac{\tilde\rho_{l,h}}{lh}$. Substituting into \eqref{pol-charge},
\bml{\ve{p}_{p,l,h}(\hat{\ve{x}}_p)&=\int_{\Box}\int_{-1}^1\ve{y}\dfrac{\tilde\rho_{l,h}}{lh}(\hat{\ve{x}}_p+l\ve{y}+hz\ve{e}_3)\dm \ve{y}\dm z=\dfrac{\tilde{\ve{p}}_{p,l,h}(\hat{\ve{x}}_p)}{lh}\\
{p}_{3,l,h}(\hat{\ve{x}}_p)&=\int_{\Box}\int_{-1}^1z\dfrac{\tilde\rho_{l,h}}{lh}(\hat{\ve{x}}_p+l\ve{y}+hz\ve{e}_3)\dm \ve{y}\dm z=\dfrac{\tilde{{p}}_{3,l,h}(\hat{\ve{x}}_p)}{lh}\\
\sigma_{l,h}(\hat{\ve{x}}_p)&=\int_{\Delta}\int_{-1}^1\dfrac{\tilde\rho_{l,h}}{lh}(\hat{\ve{x}}_p+l\ve{y}+hz\ve{e}_3)\dm \ve{y} \dm z =\dfrac{\tilde{\sigma}_{l,h}(\ve{x}_p)}{lh}\label{pol-charge-Case1}}

Substituting into \eqref{approx-sol}, we get
\bml{\Phi_{l,h}(\ve{x})&\approx\int_{\Omega}\dm\ve{x}_p' \nabla_{\ve{x}'} G(\ve{x}, {\ve{x}}'_p)\cdot\left(\tilde{\ve{p}}_{p,l,h}(\ve{x}'_p)+\dfrac{h}{l}\tilde{p}_{3,l,h}(\ve{x}'_p)\ve{e}_3\right)+\int_{\partial\Omega}\dm S_{\ve{x}_p'} G(\ve{x}, {\ve{x}}'_p)\tilde{\sigma}_{l,h}(\ve{x}'_p),\label{approx-sol-Case1}}

Setting ${h}/{l}\to0$, 
\bml{\Phi_{0}(\ve{x})&=\int_{\Omega}\dm\ve{x}_p' \nabla_{\ve{x}_p'} G(\ve{x}, {\ve{x}}'_p)\cdot \tilde{\ve{p}}_{p}(\ve{x}'_p)+\int_{\partial\Omega}\dm S_{\ve{x}_p'} G(\ve{x}, {\ve{x}}'_p)\tilde{\sigma}(\ve{x}'_p),}
where we used the fact that the polarization is in plane to replace the gradient with its tangential counterpart.

We can now perform by-parts to obtain
\bml{\Phi_{0}(\ve{x})&=\int_{\Omega}\dm\ve{x}_p'  G(\ve{x}, {\ve{x}}'_p)\divergence_p(- \tilde{\ve{p}}_{p}(\ve{x}'_p))+\int_{\partial\Omega}\dm S_{\ve{x}_p'} G(\ve{x}, {\ve{x}}'_p)\left(\tilde{\sigma}(\ve{x}'_p)+\tilde{\ve{p}}_p(\ve{x}'_p)\cdot \tilde{\ve{n}}\right)\\
&=\int_{\Omega}\dm\ve{x}_p'  G(\ve{x}, {\ve{x}}'_p)\left(-\divergence_p\tilde{\ve{p}}_{p}(\ve{x}'_p)+\mathbbm{1}_{\partial\Omega}(\tilde{\sigma}(\ve{x}'_p)+\tilde{\ve{p}}_p(\ve{x}'_p)\cdot \tilde{\ve{n}})\right),}
where $\tilde{\ve{n}}$ is tangent to $\Omega$ and normal to $\partial\Omega$.

\bml{-\Delta\Phi_{0}(\ve{x})&=-\divergence_p\tilde{\ve{p}}_p\mathbbm{1}_{\Omega}+\mathbbm{1}_{\partial\Omega}\left(\tilde{\sigma}+\tilde{\ve{p}}_p\cdot \tilde{\ve{n}}\right)~~\text{in}~\Scr{R}^3\\
\Phi_{0}(\ve{x})&=0~~\text{as}~|\ve{x}|\to\infty,\label{PDE-Case1}}

Since the potential converges weakly, the left-hand side of \eqref{PDE-Case1} is uniquely determined. Consequently, the right-hand side must also be unique in a pointwise sense. It follows that the divergence of the polarization field is uniquely determined within the bulk, whereas the surface charge density $\tilde\sigma$ is effectively compensated by the term \(\tilde{\ve{p}}_p \cdot \tilde{\ve{n}}\). Hence, it can be concluded that the non-uniqueness inherent in the tangential component of the polarization field is offset by the corresponding surface charge distribution.

\subsection{Case II. $l\ll h\to0$}\label{CaseII}
For this limit to make sense, we consider the scaling, $\rho_{l,h}=\dfrac{\tilde\rho_{l,h}}{h^2}$. Substituting into \eqref{pol-charge},
\bml{\ve{p}_{p,l,h}(\hat{\ve{x}}_p)&=\int_{\Box}\int_{-1}^1\ve{y}\dfrac{\tilde\rho_{l,h}}{h^2}(\hat{\ve{x}}_p+l\ve{y}+hz\ve{e}_3)\dm \ve{y}\dm z=\dfrac{\tilde{\ve{p}}_{p,l,h}(\hat{\ve{x}}_p)}{h^2}\\
{p}_{3,l,h}(\hat{\ve{x}}_p)&=\int_{\Box}\int_{-1}^1z\dfrac{\tilde\rho_{l,h}}{lh}(\hat{\ve{x}}_p+l\ve{y}+hz\ve{e}_3)\dm \ve{y}\dm z=\dfrac{\tilde{{p}}_{3,l,h}(\hat{\ve{x}}_p)}{h^2}\\
\sigma_{l,h}(\hat{\ve{x}}_p)&=\int_{\Delta}\int_{-1}^1\dfrac{\tilde\rho_{l,h}}{lh}(\hat{\ve{x}}_p+l\ve{y}+hz\ve{e}_3)\dm \ve{y} \dm z =\dfrac{\tilde{\sigma}_{l,h}(\ve{x}_p)}{h^2}\label{pol-charge-Case3}}

Substituting into \eqref{approx-sol}, we get
\bml{\Phi_{l,h}(\ve{x})&\approx\int_{\Omega}\dm\ve{x}_p' \nabla_{\ve{x}'} G(\ve{x}, {\ve{x}}'_p)\cdot\left(\dfrac{l}{h}\tilde{\ve{p}}_{p,l,h}(\ve{x}'_p)+\tilde{p}_{3,l,h}(\ve{x}'_p)\ve{e}_3\right)+\int_{\partial\Omega}\dm S_{\ve{x}_p'} G(\ve{x}, {\ve{x}}'_p)\dfrac{l}{h}\tilde{\sigma}_{l,h}(\ve{x}'_p),\label{approx-sol-Case3}}

Setting ${l}/{h}\to0$, 
\bml{\Phi_{0}(\ve{x})&=\int_{\Omega} \dfrac{\partial G(\ve{x}, {\ve{x}}'_p)}{\partial x_3'} \tilde{p}_{3}(\ve{x}'_p)\dm\ve{x}_p',\label{potential-Case3-i}}
which is called the electric double layer potential and is associated with the potential \cite{courant2008methods} generated by two layers of charge distribution of opposite sign and scaled magnitude brought close together. This is the dipole limit for a sheet of charge. The PDE corresponding to \eqref{potential-Case3-i} is given by  (see Appendix \ref{Appendix} for details and derivation):
\bml{\Delta\Phi_{0}(\ve{x})&=\tilde p_3(\ve{x}_p)\mathbbm{1}_{\Omega}(\ve{x}_p)\dfrac{d\delta(x_3)}{d x_3}~~\text{in}~\Scr{R}^3\\
\Phi_{0}(\ve{x})&=0~~\text{as}~|\ve{x}|\to\infty,\label{PDE-Case3}}

Given that \(\Phi_0\) is unique, it follows that \(\tilde{p}_3\) must also be uniquely determined. Consequently, the normal component of the polarization field is uniquely defined. This observation is consistent with the fact that the choice of unit cell is constrained to lie within the plane.

\subsection{Case III. $h\sim l\to0$}\label{CaseIII}
For this case, we set $h=l$ in \eqref{approx-sol}:
\bml{\Phi_{l}(\ve{x})&\approx\int_{\Omega}\dm\ve{x}_p' \nabla_{\ve{x}'} G(\ve{x}, {\ve{x}}'_p)\cdot\left(l^2\ve{p}_{p,l,h}(\ve{x}'_p)+l^2p_{3,l,h}(\ve{x}'_p)\ve{e}_3\right)+\int_{\partial\Omega}\dm S_{\ve{x}_p'} G(\ve{x}, {\ve{x}}'_p)l^2\sigma_{l,h}(\ve{x}'_p),\label{approx-sol-Case2}}

If we choose the scaling to be $\rho_l=\dfrac{\tilde{\rho}_l}{l^2}$, then from \eqref{pol-charge}, we get
\bml{\ve{p}_{p,l}(\hat{\ve{x}}_p)&=\int_{\Box}\int_{-1}^1\ve{y}\dfrac{\tilde\rho_{l,h}}{l^2}(\hat{\ve{x}}_p+l\ve{y}+hz\ve{e}_3)\dm \ve{y}\dm z=\dfrac{\tilde{\ve{p}}_{p,l}(\hat{\ve{x}}_p)}{l^2}\\
{p}_{3,l}(\hat{\ve{x}}_p)&=\int_{\Box}\int_{-1}^1z\dfrac{\tilde\rho_{l}}{l^2}(\hat{\ve{x}}_p+l\ve{y}+hz\ve{e}_3)\dm \ve{y}\dm z=\dfrac{\tilde{{p}}_{3,l}(\hat{\ve{x}}_p)}{l^2}\\
\sigma_{l}(\hat{\ve{x}}_p)&=\int_{\Delta}\int_{-1}^1\dfrac{\tilde\rho_{l}}{l^2}(\hat{\ve{x}}_p+l\ve{y}+hz\ve{e}_3)\dm \ve{y} \dm z =\dfrac{\tilde{\sigma}_{l}(\ve{x}_p)}{l^2}\label{pol-charge-Case2}}

Substituting \eqref{pol-charge-Case2} into \eqref{approx-sol-Case2} and taking the limit of $l\to0$
\bml{\Phi_{0}(\ve{x})&=\int_{\Omega}\dm\ve{x}_p' \left(\nabla_{\ve{x}'_p} G(\ve{x}, {\ve{x}}'_p)\cdot\tilde{\ve{p}}_{p}(\ve{x}'_p)+ \dfrac{\partial}{\partial x_3'}G(\ve{x}, {\ve{x}}'_p)\tilde{p}_{3}(\ve{x}'_p)\right)+\int_{\partial\Omega}\dm S_{\ve{x}_p'} G(\ve{x}, {\ve{x}}'_p)\tilde{\sigma}(\ve{x}'_p),\\
&=\int_{\Omega}\dm\ve{x}_p'  G(\ve{x}, {\ve{x}}'_p)\left(-\divergence_p\tilde{\ve{p}}_{p}(\ve{x}'_p)+\mathbbm{1}_{\partial\Omega}(\tilde{\sigma}(\ve{x}'_p)+\tilde{\ve{p}}_p(\ve{x}'_p)\cdot \tilde{\ve{n}})\right)+\dfrac{\partial}{\partial x_3'}G(\ve{x}, {\ve{x}}'_p)\tilde{p}_{3}(\ve{x}'_p)
\label{approx-sol-Case2-i}}

Therefore, the in-plane bound charge remains uniquely defined, and any non-uniqueness in the tangential component of the polarization field is compensated by the surface charge distribution. The final term, which cannot be integrated, represents the electrostatic double-layer potential and has been previously addressed in the asymptotic regime \(l \ll h \to 0\) as demonstrated in \eqref{CaseII}. Consequently, the partial differential equation, by analogy with \eqref{PDE-Case1} and \eqref{PDE-Case3}, is as follows:

\bml{-\Delta\Phi_{0}(\ve{x})&=-\divergence_p\tilde{\ve{p}}_p\mathbbm{1}_{\Omega}+\mathbbm{1}_{\partial\Omega}\left(\tilde{\sigma}+\tilde{\ve{p}}_p\cdot \tilde{\ve{n}}\right)-\tilde p_3(\ve{x}_p)\mathbbm{1}_{\Omega}(\ve{x}_p)\dfrac{d\delta(x_3)}{dx_3}~~\text{in}~\Scr{R}^3\\
\Phi_{0}(\ve{x})&=0~~\text{as}~|\ve{x}|\to\infty,\label{PDE-Case2}}

\section{Formulation}\label{PF}

For two-dimensional (2D) materials such as graphene, which consist of a single atomic layer, all atoms are confined to a 2D plane, forming a hexagonal lattice. The electronic charge distribution, however, is not restricted to this plane. For atomic configurations constrained to a 2D manifold \(\Omega \subset \Scr{R}^3\), the charge distribution can be considered localized in the vicinity of the manifold. Near the boundary of the manifold \(\partial\Omega\), the charge distribution is expected to concentrate along the corresponding curve.

One may assume that the support of the charge distribution is given by \( \Omega\oplus B^3_h \), where \( \oplus \) denotes the Minkowski sum, and \( B^3_h \) represents a ball of radius \( h \) in \( \Scr{R}^3 \). For simplicity, we approximate the support \( \Omega\oplus B^3_h \) by $(\Omega\oplus B^2_h)\times(-h,h)$. Here, \( B^2_h \) denotes a ball of radius \( h \) in \( \Scr{R}^2 \). For $h$ small, we may replace $\Omega\oplus B^2_h$ with $\tilde\Omega$ (renamed $\Omega$) and denote the support of the charge distribution as $\Omega_h=\Omega\times(-h,h)$. The advantage of this formulation is that the charge distribution is now defined over a thin film. This corresponds to the well-known thin film limit setup in mechanics \cite{james1998magnetostriction,kohn2005another,antman1995nonlinear,friesecke2002theorem}.

Extending this framework to 2D materials consisting of multiple atomic layers necessitates the identification of a representative manifold for the 2D material. This manifold may correspond to an actual atomic plane or a fictitious surface. The charge distribution is then modeled as being supported in a neighborhood of this manifold, denoted as $\Omega_h$, analogous to the previously described case. The non-uniqueness in the choice of manifold will not affect the final answer.

Due to the electrostatic attraction between positive and negative charges, materials tend to exhibit charge redistribution, with negative charge accumulating around positively charged atomic sites, thereby screening the net Coulombic force. In 3D materials, this phenomenon gives rise to a condition of charge neutrality within each unit cell. The net charge, which characterizes the electronic properties of the material, comprises both positive and negative contributions within each unit cell. As a result, the net charge distribution inherently exhibits rapid spatial variations at the microscopic scale. The dipole moment per unit cell, representing the first-order approximation to the electronic degrees of freedom, varies over macroscopic length scales, making it an ideal continuum descriptor for the electronic behavior of the material. This serves as the fundamental quantity of interest in our analysis.

For 2D materials, the unit cell is inherently two-dimensional and is defined within the plane of the chosen manifold $\Omega$, with atomic positions characterized by a motif. To define polarization, it is necessary to construct a 3D composite unit cell (normal to the manifold) for the computation of the dipole moment. Additionally, this composite unit cell must satisfy charge neutrality. For flat 2D materials, this condition is typically well satisfied due to in-plane periodicity, which ensures that charge neutrality within a composite unit cell inherently leads to charge neutrality in the bulk. For arbitrarily deformed 2D materials, defining a composite unit cell normal to the surface leads to non-uniform unit cells, complicating this approach.
However, since the parameterization space is flat, the definition of composite unit cells within this space remains straightforward.

Consider a flat 2D material that undergoes deformation from a flat initial to a curved final configuration described by the mapping $\chi:\Omega_{h,0} \to \Omega_h$. While $\Omega_{h,0}$ represents the equilibrium configuration of the electronic charge in which unit cells satisfy charge neutrality, this property does not necessarily hold for the deformed configuration $\Omega_{h}$. In response to the deformation, the electronic charge will redistribute to attain a new equilibrium state. Since Poisson’s equation is formulated with respect to the deformed configuration, our analysis begins from this configuration. Consequently, when this manifold is mapped to a planar domain, $T_h = \psi^{-1}(\Omega_h)$, the unit cells in $T$ no longer satisfy charge neutrality due to charge redistribution. Although the material as a whole remains charge neutral, individual unit cells need not maintain charge neutrality. To this end we define free charge.

\label{freecharge}
\defn{free charge of order $\alpha,\beta$}{For a charge distribution $\{\rho_{l,h}\}_{l,h}$, we define the homogenized free charge of order $\alpha,~\beta$ at a point $\ve{x}_p\in\Omega$ is defined as:
\beol{J_0(\ve{x}_p)\left(q_{l,h}\circ\psi\right)(\ve{x}_p)=\dfrac{1}{l^{\alpha}h^{\beta}}\int_{\Box}\int_{-1}^1\left(\rho_{l,h}\circ\psi\right)(\ve{x}_p+l\ve{y}+hz\ve{e}_3)\dm\ve{y},\dm z}
}
where $J_0=\sqrt{\text{det}(D\psi_0^TD\psi)}$ is the surface Jacobian for the manifold $\Omega$.

\rmrk{It is crucial to distinguish between the deformation map \( \chi \), which induces a redistribution of charge, and the parameterization map \( \psi^{-1} \), which maps the deformed configuration to \( T \) without altering the charge distribution. This distinction arises because \( \chi \) represents a physical transformation, whereas \( \psi^{-1} \) is a mathematical construct introduced for convenience. The charge redistribution resulting from \( \chi \) is governed by quantum mechanical principles and is inherently tied to the equilibrium configuration in the deformed configuration. As a consequence, the free charge density \( q_{l,h} \) in Definition \eqref{freecharge} cannot be determined from the reference configuration but must instead be derived from the deformed configuration. If $\chi$ represents an extension, the deviation from charge neutrality in the unit cells is expected to be of the order \( l^{\alpha} \). Conversely, if $\chi$ involves bending, the deviation from charge neutrality in the unit cells is anticipated to be of the order \( h^{\beta} \).}

We are now in a position to formulate the problem statement. Consider a sequence of charge distributions \( \{\rho_{l,h}\}_{l,h\in\scr{I}_1\times\scr{I}_2}\in\scr{C}(\Omega_{h},\Scr{R}) \), where \( \scr{I}_1 \) and \( \scr{I}_2 \) are indexing sets for \( l \) and \( h \), respectively. Here, \( \Omega_h = \Omega \times (-h,h) \) represents the support of the charge distribution, with \( \Omega \subset \Scr{R}^3 \) a 2D manifold. The charge distribution has free charge (Definition \ref{freecharge}) due to curvature effects.

We consider the following sequence of Poisson equations:
\bml{-\Delta\Phi_{l,h}&=\rho_{l,h}~~\text{in}~\Scr{R}^3\\
\Phi_{l,h}&=0~~\text{as}~|\ve{x}|\to\infty,\label{Problem-Curved}}

Our objective is to analyze both the continuum limit (vanishing lattice spacing) and the thin film limit of the potential induced by these charges. With two scaling parameters, we examine various asymptotic regimes of the charge density as \( l \) and \( h \) approach zero. We consider the following asymptotic regimes: \( l \ll h \to 0 \), \( l \approx h \to 0 \), and \( h \ll l \to 0 \). 

The solution to \eqref{Problem-Curved} is given by the convolution of the charge against the Green's function,
\beol{\Phi_{l,h}(\ve{r})=\int_{\Omega_h}G(\ve{r},\ve{r}')\rho_{l,h}(\ve{r}')\dm\ve{r}'\label{solution-Curved},}
where $G(\ve{r},\ve{r}'):={1}/{|\ve{r}-\ve{r}'|}$ is the free space Green's function.

Since the Green's function associated with \eqref{Problem} remains fixed, we can employ the weak convergence of \(\rho_{l,h}\) to derive a homogenized potential \(\Phi_{l,h}\)\footnote{The product of weak and strong convergence results in weak convergence.}. Our primary interest lies in analyzing the behavior of the potential in regions away from the manifold \(\Omega\). In these regions, dipoles and effective charges become evident.

\section{Homogenized Potential}\label{Homo}
In this section, we derive the homogenized potential for the three limits: \( l \ll h \to 0 \), \( l \approx h \to 0 \), and \( h \ll l \to 0 \). We divide the derivation for each into subsections. The starting point of our analysis is \eqref{solution-Curved}. For the different cases, to derive a well-defined problem, we will have to scale the charge distribution. The choice of scaling the charge distribution will differ for each asymptotic regime. 

\begin{theorem}\label{Theorem-Homo}
Let \(\{\tilde{\rho}_{l,h}\} \subset \scr{C}(\Omega_h, \Scr{R})\) denote a sequence of charge distributions, where \(\Omega\) is a two-dimensional manifold embedded in \(\Scr{R}^3\), and \(\Omega_h = \Omega \times (-h, h)\) represents the extended support of the charge distribution. Then, the following limits are satisfied:
    \begin{enumerate}
        \item Suppose $\dfrac{h}{l}\to0$ as $l,h\to0$. Consider the sequence of Poisson equations:
        \bml{-\Delta\Phi_{l,h}&=\dfrac{\rho_{l,h}}{lh}\mathbbm{1}_{\Omega_h}~~\text{in}~\Scr{R}^3\\
        \Phi_{l,h}&=0~~\text{as}~|\ve{r}|\to\infty.\label{Problem-Homogenize}}
        
        If \(\rho_{l,h}\) contains free charge of order \((1, 0)\), then the limit as \(l, h \to 0\) of \eqref{Problem-Homogenize} is given by the following partial differential equation:
        \bml{-\Delta\Phi_{0}&=\left({q}-\dfrac{\divS\left(J_0 {\ve{p}}_p\right)}{J_0}\right)\mathbbm{1}_{\Omega}+\mathbbm{1}_{\partial\Omega}\left( {\sigma}+ {\ve{p}}_p\cdot \ve{n}\right)~~\text{in}~\Scr{R}^3\\
        \Phi_{0}&=0~~\text{as}~|\ve{r}|\to\infty,\label{PDE-Case1-Homogenize}}
        where $\divS$ denotes the surface divergence operator on the manifold $\Omega$, \( J_0 \) is the Jacobian associated with surface integration over $\Omega$ (see Theorem \ref{MST} for its explicit form), $\ve{p}_p$ denotes the tangential component of the polarization field, \(q\) represents the free charge density of order \((1, 0)\), $\ve{n}$ denotes the outward-pointing unit normal vector to the boundary $\partial\Omega$ lying within the tangent plane to $\Omega$, and $\sigma$ denotes the surface charge density.

        \item Suppose $\dfrac{h}{l}\to \alpha$ finite, as $l,h\to0$. Consider the sequence of Poisson equations:
        \bml{-\Delta\Phi_{l}&=\dfrac{ \rho_{l}}{l^2}\mathbbm{1}_{\Omega_h}~~\text{in}~\Scr{R}^3\\
        \Phi_{l}&=0~~\text{as}~|\ve{r}|\to\infty.\label{Problem-Homogenize-ii}}
        
        If \(\rho_{l,h}\) possesses free charge of order \((0, 1)\) or \((1, 0)\), then the limiting behavior as \(l \to 0\) of \eqref{Problem-Homogenize-ii} is governed by the following partial differential equation:

        \bml{-\Delta\Phi_{0}&=\alpha\left({q}-\dfrac{\divS\left(J_0 {\ve{p}}_p\right)}{J_0}\right)\mathbbm{1}_{\Omega}-\alpha^2p_3\dfrac{\partial\mathbbm{1}_{\Omega}}{\partial\nu}+\mathbbm{1}_{\partial\Omega}\alpha\left( {\sigma}+ {\ve{p}}_p\cdot \ve{n}\right)~~\text{in}~\Scr{R}^3\\
        \Phi_{0}&=0~~\text{as}~|\ve{r}|\to\infty.\label{PDE-Case2-Homogenize}}
        where $\divS$ denotes the surface divergence operator on the manifold $\Omega$, \( J_0 \) is the Jacobian associated with surface integration over $\Omega$ (see Theorem \ref{MST} for its explicit form), \(\ve{p}_p\) denotes the tangential component of the polarization, \(p_3\) represents the normal component of the polarization, \(q\) refers to the free charge of order \((0,1)\) or \((1,0)\), \(\sigma\) denotes the surface charge distribution, \(\ve{n}\) is the limiting vector of tangents to \(\Omega\) normal to \(\partial\Omega\), and \(\ve{\nu}\) is the unit normal vector to the manifold \(\Omega\).

        \item Suppose $\dfrac{l}{h}\to0$ as $l,h\to0$. Consider the sequence of Poisson equations:
        \bml{-\Delta\Phi_{l}&=\dfrac{ \rho_{l}}{h^2}\mathbbm{1}_{\Omega_h}~~\text{in}~\Scr{R}^3\\
        \Phi_{l}&=0~~\text{as}~|\ve{r}|\to\infty.\label{Problem-Homogenize-iii}}
        
        If \(\rho_{l,h}\) contains a free charge of order \((0,1)\), then the limit as \(l, h \to 0\) of \eqref{Problem-Homogenize-iii} results in the following partial differential equation:
        \bml{-\Delta\Phi_{0}&= q-{p}_3\dfrac{\partial \mathbbm{1}_{\Omega}}{\partial \nu}~~\text{in}~\Scr{R}^3\\
        \Phi_{0}(\ve{r})&=0~~\text{as}~|\ve{r}|\to\infty,\label{PDE-Case3-Homogenize}}
        where \( p_3 \) denotes the normal polarization, \( q \) represents the free charge of order \((0,1)\), and \(\ve{\nu}\) denotes the unit normal vector to the manifold \(\Omega\).
        
    \end{enumerate}
\end{theorem}

In the preceding theorem, the solutions to the original PDEs converge weakly to the solutions of the corresponding homogenized limiting PDEs. Additionally, \( {\partial \mathbbm{1}_{\Omega}}/{\partial \nu} \) represents the derivative of the characteristic function in the normal direction, as elaborated in Appendix \ref{Appendix}. In what follows, we provide a rigorous proof of Theorem \ref{Theorem-Homo}, structuring the exposition such that each subpart of the theorem is addressed in a dedicated section to ensure clarity and comprehensiveness.

\subsection{Proof of Theorem \ref{Theorem-Homo}(1)}\label{Proof1}
\begin{proof}
Utilizing the Green's function, the potential \(\Phi_{l,h}\) in \eqref{Problem-Homogenize} can be expressed as:
\beol{\Phi_{l,h}(\ve{r})=\int_{\Omega}G(\ve{r},\ve{r}')\dfrac{{\rho}_{l,h}}{lh}(\ve{r}')\dm\ve{r}'\label{solution-CaseI-i}}

By transforming into the parameter space and applying Theorem \ref{MST}, equation \eqref{solution-CaseI-i} simplifies to:
\beol{\left(\Phi_{l,h}\circ\psi\right)(\ve{x})=\int_{T\times(-h,h)}G\left(\psi(\ve{x}),\psi(\ve{x}')\right)\dfrac{\left(\rho_{l,h}\circ\psi\right)(\ve{x}')}{lh}J(\ve{x}')\dm\ve{x}',\label{Transformed-Soln-CaseI}}
where $J(\ve{x})=|\text{det}D\psi(\ve{x})|=\textstyle\sqrt{\text{det}(D\psi^TD\psi)}$, $T\times(-h,h)=\psi^{-1}(\Omega_h)=:\psi^{-1}\left(\Omega\times(-h,h)\right)$.

For brevity, we define the following composite variables; $\tilde\Phi_{l,h}=\left(\Phi_{l,h}\circ\psi\right)$, $\tilde{G}=G\circ\psi$ and $\tilde{\rho}_{l,h}=\left(\rho_{l,h}\circ\psi\right) $. With respect to these composite variables \eqref{Transformed-Soln-CaseI} takes the form
\beol{\tilde{\Phi}_{l,h}(\ve{x})=\int_{T_h}\tilde{G}(\ve{x},\ve{x}')\dfrac{\tilde{\rho}_{l,h}(\ve{x}')}{lh}J(\ve{x}')\dm\ve{x}'\label{Transformed-Soln-CaseI-ii}}

We now decompose the integral in \eqref{Transformed-Soln-CaseI-ii} over the domain $T$ into the disjoint regions in $T(\Box)\cup T(\Delta)$. This is achieved by applying the corner map substitution:  
\beol{
\ve{x}=\hat{\ve{x}}_p+l\ve{y}+hz\ve{e}_3,\label{corner-map-sub-CaseI}}
where $\ve{x}_p\in T$ denotes the projection of $\ve{x}$ onto $T$. The normal direction $\ve{e}_3$ is scaled by $h$, yielding  $\ve{x}=\ve{x}_p+x_3\ve{e}_3=\ve{x}_p+hz\ve{e}_3$. The corner map $\hat{\ve{x}}_p$ (see Section \ref{Notation}) enables us to express $\ve{x}_p=\hat{\ve{x}}_p+l\ve{y}$ with $\ve{y}\in\Box$.

Using \eqref{corner-map-sub-CaseI},
\bml{\tilde{\Phi}_{l,h}(\ve{x})&=\sum_{\hat{\ve{x}}'_p\in T(\Box)}l^2\int_{\Box}\int_{-1}^1\tilde{G}(\ve{x},\hat{\ve{x}}'_p+l\ve{y}'+hz'\ve{e}_3)\dfrac{\tilde{\rho}_{l,h}(\hat{\ve{x}}'_p+l\ve{y}'+hz'\ve{e}_3)}{l}J(\hat{\ve{x}}'_p+l\ve{y}'+hz'\ve{e}_3)\dm\ve{y}'\dm z'\\
&+\sum_{\hat{\ve{x}}'_p\in T(\Delta)}l\int_{\Box}\int_{-1}^1\tilde{G}(\ve{x},\hat{\ve{x}}'_p+l\ve{y}'+hz'\ve{e}_3){\tilde{\rho}_{l,h}(\hat{\ve{x}}'_p+l\ve{y}'+hz'\ve{e}_3)}(\hat{\ve{x}}'_p+l\ve{y}'+hz'\ve{e}_3)\dm\ve{y}'\dm z'\\
&=\text{term I} + \text{term II}\label{CornerMap-Integration-CaseI}}

We aim to homogenize the potential away from $T$. Since the Green's function is smooth away from $\Omega$, it admits a Taylor series expansion. Using the condition of total derivative, we can rewrite the first term in \eqref{CornerMap-Integration-CaseI},
\bml{\text{term I}&=\sum_{\hat{\ve{x}}'_p\in T(\Box)}l^2\int_{\Box}\int_{-1}^1\left[\dfrac{\tilde{G}(\ve{x},\ve{x}')-\tilde{G}(\ve{x},\hat{\ve{x}}'_p)-\left(l\ve{y}'+hz'\ve{e}_3\right)\cdot \nabla_{\ve{x}'}\tilde{G}(\ve{x},\hat{\ve{x}}'_p)}{|l\ve{y}'+hz'\ve{e}_3|}\right]|\ve{y}'+\dfrac{h}{l}z'\ve{e}_3|\tilde{\rho}_{l,h}(\ve{x}')J(\ve{x}')\dm\ve{y}'\dm z'\\
&+\sum_{\hat{\ve{x}}'_p\in T(\Box)}l^2\tilde{G}(\ve{x},\hat{\ve{x}}'_p)\int_{\Box}\int_{-1}^1\dfrac{\tilde{\rho}_{l,h}(\ve{x}')}{l}J(\ve{x}')\dm\ve{y}'\dm z'\\
&+\sum_{\hat{\ve{x}}'_p\in T(\Box)}l^2\nabla_{\ve{x}'}\tilde{G}(\ve{x},\hat{\ve{x}}'_p)\cdot\int_{\Box}\int_{-1}^1\left(\ve{y}'+\dfrac{h}{l}z'\ve{e}_3\right)\tilde{\rho}_{l,h}(\ve{x}')J(\ve{x}')\dm\ve{y}'\dm z'.\label{Taylor-Integration-Bulk-CaseI}}

The second term is a free charge of order $1,0$ (see definition \ref{freecharge}). The third term gives us polarization. We rewrite the above as
\bml{\text{term I}&=\sum_{\hat{\ve{x}}'_p\in T(\Box)}l^2\int_{\Box}\int_{-1}^1\left[\dfrac{\tilde{G}(\ve{x},\ve{x}')-\tilde{G}(\ve{x},\hat{\ve{x}}'_p)-\left(l\ve{y}'+hz'\ve{e}_3\right)\cdot \nabla_{\ve{x}'}\tilde{G}(\ve{x},\hat{\ve{x}}'_p)}{|l\ve{y}'+hz'\ve{e}_3|}\right]|\ve{y}'+\dfrac{h}{l}z'\ve{e}_3|\tilde{\rho}_{l,h}(\ve{x}')J(\ve{x}')\dm\ve{y}'\dm z'\\
&+\sum_{\hat{\ve{x}}'_p\in T(\Box)}l^2\tilde{G}(\ve{x},\hat{\ve{x}}'_p)\tilde{q}_{l,h}(\hat{\ve{x}}'_p)J_0(\hat{\ve{x}}'_p)+\sum_{\hat{\ve{x}}'_{p}\in T(\Box)}l^2\nabla_{\ve{x}'}\tilde{G}(\ve{x},\hat{\ve{x}}'_p)\cdot\left(\tilde{\ve{p}}_{p,l,h}(\hat{\ve{x}}'_p)+\dfrac{h}{l}\tilde p_{3,l,h}(\hat{\ve{x}}'_p)\ve{e}_3\right)J_0(\hat{\ve{x}}'_p),\label{Taylor-Integration-Bulk-CaseI-ii}}
where we defined
\bml{J_0(\hat{\ve{x}}'_p)\tilde{\ve{p}}_{p,l,h}(\hat{\ve{x}}'_p)&=\int_{\Box}\int_{-1}^1 \ve{y}'\tilde{\rho}_{l,h}(\hat{\ve{x}}'_p+l\ve{y}'+hz'\ve{e}_3)J(\hat{\ve{x}}'_p+l\ve{y}'+hz'\ve{e}_3)\dm\ve{y}'\dm z'\\
J_0(\hat{\ve{x}}'_p)\tilde p_{3,l,h}(\hat{\ve{x}}'_p)&=\int_{\Box}\int_{-1}^1z'\tilde{\rho}_{l,h}(\hat{\ve{x}}'_p+l\ve{y}'+hz'\ve{e}_3)J(\hat{\ve{x}}'_p+l\ve{y}'+hz'\ve{e}_3)\dm\ve{y}'\dm z',\\
}
where $J_0=\sqrt{\text{det}(D\psi_0^TD\psi_0)}$ is the Jacobian for the integration on the manifold $\Omega$. Since $\ve{r}\not\in\Omega$, $G(\ve{r},\ve{r}')$ does not contain a singularity. Thus, $G$ has the required smoothness for us to compute derivatives. Noting that ${h}/{l}\to0$ in this limit, the term $\ve{y}'+({h}/{l})z'\ve{e}_3\to \ve{y}' $ a finite limit. Since $\tilde{G}$ is differentiable, the first term converges to zero from the definition of the derivative. Further, note that $\tilde{\rho}_{l,h}$ and $J$ are bounded and $\textstyle\sum_{T(\Box)} l^2\to\textstyle\int_T\dm\ve{x}'_p$. Thus, the limit of the first term is zero. Taking the limit $h,l\to0$, we get
\beol{\text{term I}\mapsto\int_T\left(\tilde{G}(\ve{x},\ve{x}'_p)\tilde{q}(\ve{x}'_p)+\nabla_{\ve{x}'_p}\tilde{G}(\ve{x},\ve{x}'_p)\cdot\tilde{\ve{p}}_p(\ve{x}'_p)\right)J_0(\hat{\ve{x}}'_p)\dm\ve{x}'_p,\label{BulkLimit-CaseI}}
where in the limit the normal polarization $p_3$ disappears and $\hat{\ve{x}}_p\mapsto\ve{x}_p$  giving us an integral over $T$.

Similarly, using the conditions of continuity, we can rewrite the second term in \eqref{CornerMap-Integration-CaseI} as
\bml{\text{term II}&=\sum_{\hat{\ve{x}}'_p\in T(\Delta)}l \int_{\Delta}\int_{-1}^1\left[\tilde{G}(\ve{x},\hat{\ve{x}}'_p+l\ve{y}'+hz'\ve{e}_3)-\tilde{G}(\ve{x},\hat{\ve{x}}'_p)\right]\tilde{\rho}_{l,h}({\ve{x}}')J({\ve{x}}')\dm\ve{y}'\dm z'\\
&+\sum_{\hat{\ve{x}}'_p\in T(\Delta)}l \tilde{G}(\ve{x},\hat{\ve{x}}'_p)\int_{\Delta}\int_{-1}^1\tilde{\rho}_{l,h}(\hat{\ve{x}}'_p+l\ve{y}'+hz'\ve{e}_3)J(\hat{\ve{x}}'_p+l\ve{y}'+hz'\ve{e}_3)\dm\ve{y}'\dm z'.\label{Taylor-Integration-Surface-CaseI}}

From continuity, the first term in \eqref{Taylor-Integration-Surface-CaseI} vanishes in the limit. The $\textstyle\sum_{T(\Delta)}l\mapsto\int_{\partial T}\dm S_{\ve{x}'_p}$. 

We define
\beol{J_0(\hat{\ve{x}}'_p)\tilde{\sigma}(\hat{\ve{x}}'_p)=\int_{\Delta}\int_{-1}^1\tilde{\rho}_{l,h}(\hat{\ve{x}}'_p+l\ve{y}'+hz'\ve{e}_3)J(\hat{\ve{x}}'_p+l\ve{y}'+hz'\ve{e}_3)\dm\ve{y}'\dm z'.}

Taking the limit of $l,h\to0$, we get
\beol{\text{term II}\mapsto\int_{\partial T}\tilde{G}(\ve{x},\ve{x}'_p)\tilde{\sigma}(\ve{x}'_p)J_0(\ve{x}'_p)\dm S_{\ve{x}'_p}.\label{SurfaceLimit-CaseI}}

Putting \eqref{BulkLimit-CaseI} and \eqref{SurfaceLimit-CaseI} together, we get
\bml{\tilde{\Phi}_{0}(\ve{x})&=\int_{T}\left(\tilde{G}(\ve{x},\ve{x}'_p)\tilde{q}(\ve{x}'_p)+\nabla_{\ve{x}'_p}\tilde{G}(\ve{x},\ve{x}'_p)\cdot\tilde{\ve{p}}_p(\ve{x}'_p)\right)J_0( {\ve{x}}'_p)\dm\ve{x}'_p+\int_{\partial T}\tilde{G}(\ve{x},\ve{x}'_p)\tilde{\sigma}(\ve{x}'_p)J_0(\ve{x}'_p)\dm S_{\ve{x}'_p}\\
&=\int_{T}\tilde{G}(\ve{x},\ve{x}'_p)\left(\tilde{q}(\ve{x}'_p)-\dfrac{\divergence_p\left(J_0\tilde{\ve{p}}_p(\ve{x}'_p)\right)}{J_0( {\ve{x}}'_p)}\right)J_0( {\ve{x}}'_p)\dm\ve{x}'_p+\int_{\partial T}\tilde{G}(\ve{x},\ve{x}'_p)\left(\tilde{\sigma}(\ve{x}'_p)+\tilde{\ve{p}}_p(\ve{x}'_p)\cdot \tilde{\ve{n}}\right)J_0(\ve{x}'_p)\dm S_{\ve{x}'_p}\\
&=\int_{T}\tilde{G}(\ve{x},\ve{x}'_p)\left(\tilde{q}(\ve{x}'_p)-\dfrac{\divergence_p\left(J_0\tilde{\ve{p}}_p(\ve{x}'_p)\right)}{J_0( {\ve{x}}'_p)}+\mathbbm{1}_{\partial T}\left(\tilde{\sigma}(\ve{x}'_p)+\tilde{\ve{p}}_p(\ve{x}'_p)\cdot \tilde{\ve{n}}\right)\right)J_0( {\ve{x}}'_p)\dm\ve{x}'_p,}
where $\tilde{\ve{n}}$ is the normal to $\partial T$.

We can now transform the integral back into real space to get
\beol{\Phi(\ve{r})=\int_{\Omega} {G}(\ve{r},\ve{r}')\left( {q}(\ve{r}')-\dfrac{\divS\left(J_0 {\ve{p}}_p\right)}{J_0}+\mathbbm{1}_{\partial \Omega}\left( {\sigma}(\ve{r}')+ {\ve{p}}_p(\ve{r}')\cdot \ve{n}\right)\right)\dm\ve{r}',\label{Limit-Integral-CaseI}}
where $\divS$ is the surface divergence and $\ve{n}$ is the limit of tangents to $\Omega$ normal to $\partial\Omega$.

The above is the solution to the following PDE
\bml{-\Delta\Phi_{0}&=\left({q}-\dfrac{\divS\left(J_0 {\ve{p}}_p\right)}{J_0}\right)\mathbbm{1}_{\Omega}+\mathbbm{1}_{\partial\Omega}\left( {\sigma}+ {\ve{p}}_p\cdot \ve{n}\right)~~\text{in}~\Scr{R}^3\\
        \Phi_{0}&=0~~\text{as}~|\ve{r}|\to\infty.\label{PDE-Case1-Homogenized}}
which is the result in \eqref{PDE-Case1-Homogenize}.
\end{proof}

\subsection{Proof of Theorem \ref{Theorem-Homo}(2)}\label{Proof2}
\begin{proof}
Making use of the Green's function, we can write $\Phi_{l,h}$ in \eqref{Problem-Homogenize-ii} as:
\beol{\Phi_{l}(\ve{r})=\int_{\Omega}G(\ve{r},\ve{r}')\dfrac{{\rho}_{l}}{l^2}(\ve{r}')\dm\ve{r}'\label{solution-CaseII-i}}

Transforming into the parameter space and utilizing theorem \ref{MST}, we simplify equation \eqref{solution-CaseII-i} to obtain:
\beol{\left(\Phi_{l}\circ\psi\right)(\ve{x})=\int_{T\times(-h,h)}G\left(\psi(\ve{x}),\psi(\ve{x}')\right)\dfrac{\left(\rho_{l}\circ\psi\right)(\ve{x}')}{l^2}J(\ve{x}')\dm\ve{x}',\label{Transformed-Soln-CaseII}}
where $J(\ve{x})=|\text{det}D\psi(\ve{x})|$, $T\times(-h,h)=\psi^{-1}(\Omega_h)=:\psi^{-1}\left(\Omega\times(-h,h)\right)$. Since we are interested in $h/l\to\alpha$, we can replace $h$ with $\alpha l$.

For brevity, we define the following composite variables; $\tilde\Phi_{l}=\left(\Phi_{l}\circ\psi\right)$, $\tilde{G}=G\circ\psi$ and $\tilde{\rho}_{l}=\left(\rho_{l}\circ\psi\right) $. With respect to these composite variables \eqref{Transformed-Soln-CaseII} takes the form
\beol{\tilde{\Phi}_{l}(\ve{x})=\int_{T\times(-\alpha l,\alpha l)}\tilde{G}(\ve{x},\ve{x}')\dfrac{\tilde{\rho}_{l}(\ve{x}')}{l^2}J(\ve{x}')\dm\ve{x}'\label{Transformed-Soln-CaseII-ii}}

We now decompose the integral in \eqref{Transformed-Soln-CaseII-ii} over the domain $T$ into the disjoint regions in $T(\Box)\cup T(\Delta)$. This is achieved by applying the corner map substitution:  
\beol{
\ve{x}=\hat{\ve{x}}_p+l\ve{y}+\alpha lz\ve{e}_3,\label{corner-map-sub-CaseII}}
where $\ve{x}_p\in T$ denotes the projection of $\ve{x}$ onto $T$. The normal direction $\ve{e}_3$ is scaled by $h=\alpha l$, yielding  $\ve{x}=\ve{x}_p+x_3\ve{e}_3=\ve{x}_p+\alpha lz\ve{e}_3$. The corner map $\hat{\ve{x}}_p$ (see Section \ref{Notation}) enables us to express $\ve{x}_p=\hat{\ve{x}}_p+l\ve{y}$ with $\ve{y}\in\Box$.

Using \eqref{corner-map-sub-CaseII},
\bml{\tilde{\Phi}_{l}(\ve{x})&=\sum_{\hat{\ve{x}}'_p\in T(\Box)}l^2\int_{\Box}\int_{-1}^1\tilde{G}(\ve{x},\hat{\ve{x}}'_p+l\ve{y}'+\alpha 
 lz'\ve{e}_3)\alpha \dfrac{\tilde{\rho}_{l}(\hat{\ve{x}}'_p+l\ve{y}'+\alpha lz'\ve{e}_3)}{l}J(\hat{\ve{x}}'_p+l\ve{y}'+\alpha 
 lz'\ve{e}_3)\dm\ve{y}'\dm z'\\
&+\sum_{\hat{\ve{x}}'_p\in T(\Delta)}l\int_{\Box}\int_{-1}^1\tilde{G}(\ve{x},\hat{\ve{x}}'_p+l\ve{y}'+\alpha 
 lz'\ve{e}_3) \alpha 
 {\tilde{\rho}_{l}(\hat{\ve{x}}'_p+l\ve{y}'+\alpha 
 lz'\ve{e}_3)}(\hat{\ve{x}}'_p+l\ve{y}'+\alpha 
 lz'\ve{e}_3)\dm\ve{y}'\dm z'\\
&=\text{term I} + \text{term II}\label{CornerMap-Integration-CaseII}}

We aim to homogenize the potential away from $T$. Since the Green's function is smooth away from $\Omega$, it admits a Taylor series expansion. Using the condition of total derivative, we can rewrite the first term in \eqref{CornerMap-Integration-CaseII},
\bml{\text{term I}&=\sum_{\hat{\ve{x}}'_p\in T(\Box)}l^2\int_{\Box}\int_{-1}^1\left[\dfrac{\tilde{G}(\ve{x},\ve{x}')-\tilde{G}(\ve{x},\hat{\ve{x}}'_p)-\left(l\ve{y}'+\alpha lz'\ve{e}_3\right)\cdot \nabla_{\ve{x}'}\tilde{G}(\ve{x},\hat{\ve{x}}'_p)}{|l\ve{y}'+\alpha lz'\ve{e}_3|}\right]|\ve{y}'+\alpha z'\ve{e}_3|\alpha\tilde{\rho}_{l}(\ve{x}')J(\ve{x}')\dm\ve{y}'\dm z'\\
&+\sum_{\hat{\ve{x}}'_p\in T(\Box)}l^2\tilde{G}(\ve{x},\hat{\ve{x}}'_p)\int_{\Box}\int_{-1}^1\alpha\dfrac{\tilde{\rho}_{l}(\ve{x}')}{l}J(\ve{x}')\dm\ve{y}'\dm z'\\
&+\sum_{\hat{\ve{x}}'_p\in T(\Box)}l^2h\nabla_{\ve{x}'}\tilde{G}(\ve{x},\hat{\ve{x}}'_p)\cdot\int_{\Box}\int_{-1}^1\left(\ve{y}'+\alpha z'\ve{e}_3\right)\alpha\tilde{\rho}_{l}(\ve{x}')J(\ve{x}')\dm\ve{y}'\dm z'.\label{Taylor-Integration-Bulk-CaseII}}

The second term is a free charge of order $(1,0)$ or $(0,1)$ (see definition \ref{freecharge}). The third term gives us polarization. We rewrite the above as
\bml{\text{term I}&=\sum_{\hat{\ve{x}}'_p\in T(\Box)}l^2\int_{\Box}\int_{-1}^1\left[\dfrac{\tilde{G}(\ve{x},\ve{x}')-\tilde{G}(\ve{x},\hat{\ve{x}}'_p)-\left(l\ve{y}'+\alpha lz'\ve{e}_3\right)\cdot \nabla_{\ve{x}'}\tilde{G}(\ve{x},\hat{\ve{x}}'_p)}{|l\ve{y}'+\alpha lz'\ve{e}_3|}\right]|\ve{y}'+\alpha z'\ve{e}_3|\alpha\tilde{\rho}_{l}(\ve{x}')J(\ve{x}')\dm\ve{y}'\dm z'\\
&+\sum_{\hat{\ve{x}}'_p\in T(\Box)}l^2\tilde{G}(\ve{x},\hat{\ve{x}}'_p)\alpha\tilde{q}_{l}(\hat{\ve{x}}'_p)J_0(\hat{\ve{x}}'_p)+\sum_{\hat{\ve{x}}'_{p}\in T(\Box)}l^2\nabla_{\ve{x}'}\tilde{G}(\ve{x},\hat{\ve{x}}'_p)\cdot\left(\alpha\tilde{\ve{p}}_{p,l}(\hat{\ve{x}}'_p)+\alpha^2\tilde p_{3,l}(\hat{\ve{x}}'_p)\ve{e}_3\right)J_0(\hat{\ve{x}}'_p),\label{Taylor-Integration-Bulk-CaseII-ii}}
where we defined
\bml{J_0(\hat{\ve{x}}'_p)\tilde{\ve{p}}_{p,l}(\hat{\ve{x}}'_p)&=\int_{\Box}\int_{-1}^1 \ve{y}'\tilde{\rho}_{l}(\hat{\ve{x}}'_p+l\ve{y}'+\alpha lz'\ve{e}_3)J(\hat{\ve{x}}'_p+l\ve{y}'+\alpha lz'\ve{e}_3)\dm\ve{y}'\dm z',\\
J_0(\hat{\ve{x}}'_p)\tilde p_{3,l}(\hat{\ve{x}}'_p)&=\int_{\Box}\int_{-1}^1z'\tilde{\rho}_{l}(\hat{\ve{x}}'_p+l\ve{y}'+\alpha lz'\ve{e}_3)J(\hat{\ve{x}}'_p+l\ve{y}'+\alpha lz'\ve{e}_3)\dm\ve{y}'\dm z',}
where $J_0=\textstyle\sqrt{\text{det}(D\psi_0^TD\psi_0)}$ is the Jacobian associated with integration on the manifold $\Omega$. Since $\ve{r}\not\in\Omega$, $G(\ve{r},\ve{r}')$ does not contain a singularity. Thus, $G$ has the required smoothness for us to compute derivatives. Since $\tilde{G}$ is differentiable, the first term converges to zero from the definition of the derivative. Further, note that $\tilde{\rho}_{l}$ and $J$ are bounded and $\textstyle\sum_{T(\Box)} l^2\to\textstyle\int_T\dm\ve{x}'_p$. Thus, the limit of the first term is zero. Taking the limit $l\to0$, we get
\beol{\text{term I}\mapsto\int_T\left(\tilde{G}(\ve{x},\ve{x}'_p)\alpha\tilde{q}(\ve{x}'_p)+\nabla_{\ve{x}'}\tilde{G}(\ve{x},\ve{x}'_p)\alpha\cdot\tilde{\ve{p}}_p(\ve{x}'_p)+\alpha^2\dfrac{\partial}{\partial x_3'}\tilde{G}(\ve{x},\ve{x}'_p)\tilde{p}_3(\ve{x}'_p)\right)J_0(\hat{\ve{x}}'_p)\dm\ve{x}'_p.\label{BulkLimit-CaseII}}

Similarly, using the conditions of continuity, we can rewrite the second term in \eqref{CornerMap-Integration-CaseII} as
\bml{\text{term II}&=\sum_{\hat{\ve{x}}'_p\in T(\Delta)}l \int_{\Delta}\int_{-1}^1\left[\tilde{G}(\ve{x},\hat{\ve{x}}'_p+l\ve{y}'+\alpha lz'\ve{e}_3)-\tilde{G}(\ve{x},\hat{\ve{x}}'_p)\right]\alpha\tilde{\rho}_{l}({\ve{x}}')J({\ve{x}}')\dm\ve{y}'\dm z'\\
&+\sum_{\hat{\ve{x}}'_p\in T(\Delta)}l \tilde{G}(\ve{x},\hat{\ve{x}}'_p)\int_{\Delta}\int_{-1}^1\alpha\tilde{\rho}_{l}(\hat{\ve{x}}'_p+l\ve{y}'+\alpha lz'\ve{e}_3)J(\hat{\ve{x}}'_p+l\ve{y}'+\alpha lz'\ve{e}_3)\dm\ve{y}'\dm z'.\label{Taylor-Integration-Surface-CaseII}}

From continuity, the first term in \eqref{Taylor-Integration-Surface-CaseII} vanishes in the limit. The $\textstyle\sum_{T(\Delta)}l\mapsto\int_{\partial T}\dm S_{\ve{x}'_p}$. We define
\beol{J_0(\hat{\ve{x}}'_p)\tilde{\sigma}(\hat{\ve{x}}'_p)=\int_{\Delta}\int_{-1}^1\tilde{\rho}_{l}(\hat{\ve{x}}'_p+l\ve{y}'+\alpha lz'\ve{e}_3)J(\hat{\ve{x}}'_p+l\ve{y}'+\alpha lz'\ve{e}_3)\dm\ve{y}'\dm z'.}

Taking the limit of $l\to0$, we get
\beol{\text{term II}\mapsto\int_{\partial T}\tilde{G}(\ve{x},\ve{x}'_p)\alpha \tilde{\sigma}(\ve{x}'_p)J_0(\ve{x}'_p)\dm S_{\ve{x}'_p}.\label{SurfaceLimit-CaseII}}

Putting \eqref{BulkLimit-CaseII} and \eqref{SurfaceLimit-CaseII} together, we get
\bml{\tilde{\Phi}_{0}(\ve{x})&=\int_{T}\left(\tilde{G}(\ve{x},\ve{x}'_p)\alpha\tilde{q}(\ve{x}'_p)+\nabla_{\ve{x}'_p}\tilde{G}(\ve{x},\ve{x}'_p)\cdot\alpha\tilde{\ve{p}}_p(\ve{x}'_p)+\alpha^2\dfrac{\partial}{\partial x_3'}\tilde{G}(\ve{x},\ve{x}'_p)\tilde{p}_3(\ve{x}'_p)\right)J_0( {\ve{x}}'_p)\dm\ve{x}'_p\\
&+\int_{\partial T}\tilde{G}(\ve{x},\ve{x}'_p)\alpha\tilde{\sigma}(\ve{x}'_p)J_0(\ve{x}'_p)\dm S_{\ve{x}'_p}\\
&=\int_{T}\tilde{G}(\ve{x},\ve{x}'_p)\left(\alpha\tilde{q}(\ve{x}'_p)-\dfrac{\divergence_p\left(J_0\alpha\tilde{\ve{p}}_p(\ve{x}'_p)\right)}{J_0( {\ve{x}}'_p)}\right)J_0( {\ve{x}}'_p)+\alpha^2\dfrac{\partial\tilde{G}(\ve{x},\ve{x}'_p)}{\partial x_3'}\tilde{p}_3(\ve{x}'_p)J_0( {\ve{x}}'_p)\dm\ve{x}'_p\\
&+\int_{\partial T}\tilde{G}(\ve{x},\ve{x}'_p)\alpha\left(\tilde{\sigma}(\ve{x}'_p)+\tilde{\ve{p}}_p(\ve{x}'_p)\cdot \tilde{\ve{n}}\right)J_0(\ve{x}'_p)\dm S_{\ve{x}'_p}\\
&=\int_{T}\tilde{G}(\ve{x},\ve{x}'_p)\left(\alpha\tilde{q}-\alpha\dfrac{\divergence_p\left(J_0\tilde{\ve{p}}_p\right)}{J_0}+\mathbbm{1}_{\partial T}\alpha\left(\tilde{\sigma}+\tilde{\ve{p}}_p\cdot \tilde{\ve{n}}\right)\right)J_0( {\ve{x}}'_p)+\alpha^2\dfrac{\partial\tilde{G}(\ve{x},\ve{x}'_p)}{\partial x_3'}\tilde{p}_3(\ve{x}'_p)J_0(\ve{x}'_p)\dm\ve{x}'_p,}
where $\tilde{\ve{n}}$ is normal to $\partial T$.

We can now transform the integral back into real space to get
\beol{\Phi(\ve{r})=\int_{\Omega} {G}(\ve{r},\ve{r}')\left(\alpha{q}(\ve{r}')-\alpha\dfrac{\divS\left(J_0 {\ve{p}}_p\right)}{J_0}+\mathbbm{1}_{\partial \Omega}\alpha\left( {\sigma}(\ve{r}')+ {\ve{p}}_p(\ve{r}')\cdot \ve{n}\right)\right)+\alpha^2\dfrac{\partial {G}(\ve{r},\ve{r}')}{\partial \nu'}\tilde{p}_3(\ve{r}')\dm\ve{r}',\label{Limit-Integral-CaseII}}
where $\divS$ is the surface divergence and $\ve{n}$ is the limit of tangents to $\Omega$ normal to $\partial\Omega$.

The above is the solution to the following PDE
\bml{-\Delta\Phi_{0}&=\alpha\left({q}-\dfrac{\divS\left(J_0 {\ve{p}}_p\right)}{J_0}\right)\mathbbm{1}_{\Omega}-\alpha^2p_3\dfrac{\partial\mathbbm{1}_{\Omega}}{\partial\nu}+\mathbbm{1}_{\partial\Omega}\alpha\left( {\sigma}+ {\ve{p}}_p\cdot \ve{n}\right)~~\text{in}~\Scr{R}^3\\
        \Phi_{0}&=0~~\text{as}~|\ve{r}|\to\infty.\label{PDE-Case2-Homogenized}}
which is the result in \eqref{PDE-Case2-Homogenize}.
\end{proof}

\subsection{Proof of Theorem \ref{Theorem-Homo}(3)}\label{Proof3}
\begin{proof}
Making use of the Green's function, we can write $\Phi_{l,h}$ in \eqref{Problem-Homogenize-iii} as:
\beol{\Phi_{l,h}(\ve{r})=\int_{\Omega}G(\ve{r},\ve{r}')\dfrac{{\rho}_{l,h}}{h^2}(\ve{r}')\dm\ve{r}'\label{solution-CaseIII-i}}

Transforming into the parameter space and utilizing theorem \ref{MST}, we simplify equation \eqref{solution-CaseIII-i} to obtain:
\beol{\left(\Phi_{l,h}\circ\psi\right)(\ve{x})=\int_{T\times(-h,h)}G\left(\psi(\ve{x}),\psi(\ve{x}')\right)\dfrac{\left(\rho_{l,h}\circ\psi\right)(\ve{x}')}{h^2}J(\ve{x}')\dm\ve{x}',\label{Transformed-Soln-CaseIII}}
where $J(\ve{x})=|\text{det}D\psi(\ve{x})|$, $T\times(-h,h)=\psi^{-1}(\Omega_h)=:\psi^{-1}\left(\Omega\times(-h,h)\right)$.

For brevity, we define the following composite variables; $\tilde\Phi_{l,h}=\left(\Phi_{l,h}\circ\psi\right)$, $\tilde{G}=G\circ\psi$ and $\tilde{\rho}_{l,h}=\left(\rho_{l,h}\circ\psi\right) $. With respect to these composite variables \eqref{Transformed-Soln-CaseIII} takes the form
\beol{\tilde{\Phi}_{l,h}(\ve{x})=\int_{T\times(-h,h)}\tilde{G}(\ve{x},\ve{x}')\dfrac{\tilde{\rho}_{l,h}(\ve{x}')}{h^2}J(\ve{x}')\dm\ve{x}'\label{Transformed-Soln-CaseIII-ii}}

We now decompose the integral in \eqref{Transformed-Soln-CaseIII-ii} over the domain $T$ into the disjoint regions in $T(\Box)\cup T(\Delta)$. This is achieved by applying the corner map substitution:  
\beol{
\ve{x}=\hat{\ve{x}}_p+l\ve{y}+hz\ve{e}_3,\label{corner-map-sub-CaseIII}}
where $\ve{x}_p\in T$ denotes the projection of $\ve{x}$ onto $T$. The normal direction $\ve{e}_3$ is scaled by $h$, yielding  $\ve{x}=\ve{x}_p+x_3\ve{e}_3=\ve{x}_p+hz\ve{e}_3$. The corner map $\hat{\ve{x}}_p$ (see Section \ref{Notation}) enables us to express $\ve{x}_p=\hat{\ve{x}}_p+l\ve{y}$ with $\ve{y}\in\Box$.

Using \eqref{corner-map-sub-CaseIII},
\bml{\tilde{\Phi}_{l,h}(\ve{x})&=\sum_{\hat{\ve{x}}'_p\in T(\Box)}l^2\int_{\Box}\int_{-1}^1\tilde{G}(\ve{x},\hat{\ve{x}}'_p+l\ve{y}'+hz'\ve{e}_3)\dfrac{\tilde{\rho}_{l,h}(\hat{\ve{x}}'_p+l\ve{y}'+hz'\ve{e}_3)}{h}J(\hat{\ve{x}}'_p+l\ve{y}'+hz'\ve{e}_3)\dm\ve{y}'\dm z'\\
&+\sum_{\hat{\ve{x}}'_p\in T(\Delta)}l\dfrac{l}{h}\int_{\Box}\int_{-1}^1\tilde{G}(\ve{x},\hat{\ve{x}}'_p+l\ve{y}'+hz'\ve{e}_3){\tilde{\rho}_{l,h}(\hat{\ve{x}}'_p+l\ve{y}'+hz'\ve{e}_3)}(\hat{\ve{x}}'_p+l\ve{y}'+hz'\ve{e}_3)\dm\ve{y}'\dm z'\\
&=\text{term I} + \text{term II}\label{CornerMap-Integration-CaseIII}}

Observe that the sum over \( T(\Delta) \) is of order \({l}/{h}\), which vanishes in the limit as \( l \to 0 \). Consequently, our focus shifts to term I. The objective is to homogenize the potential away from \( T \). Since the Green's function remains smooth in regions away from \( \Omega \), it admits a Taylor series expansion. By leveraging the condition of total differentiability, we can reformulate the first term in \eqref{CornerMap-Integration-CaseIII},
\bml{\text{term I}&=\sum_{\hat{\ve{x}}'_p\in T(\Box)}l^2\int_{\Box}\int_{-1}^1\left[\dfrac{\tilde{G}(\ve{x},\ve{x}')-\tilde{G}(\ve{x},\hat{\ve{x}}'_p)-\left(l\ve{y}'+hz'\ve{e}_3\right)\cdot \nabla_{\ve{x}'}\tilde{G}(\ve{x},\hat{\ve{x}}'_p)}{|l\ve{y}'+hz'\ve{e}_3|}\right]|\dfrac{l}{h}\ve{y}'+z'\ve{e}_3|\tilde{\rho}_{l,h}(\ve{x}')J(\ve{x}')\dm\ve{y}'\dm z'\\
&+\sum_{\hat{\ve{x}}'_p\in T(\Box)}l^2\tilde{G}(\ve{x},\hat{\ve{x}}'_p)\int_{\Box}\int_{-1}^1\dfrac{\tilde{\rho}_{l,h}(\ve{x}')}{h}J(\ve{x}')\dm\ve{y}'\dm z'\\
&+\sum_{\hat{\ve{x}}'_p\in T(\Box)}l^2\nabla_{\ve{x}'}\tilde{G}(\ve{x},\hat{\ve{x}}'_p)\cdot\int_{\Box}\int_{-1}^1\left(\dfrac{l}{h}\ve{y}'+z'\ve{e}_3\right)\tilde{\rho}_{l,h}(\ve{x}')J(\ve{x}')\dm\ve{y}'\dm z'.\label{Taylor-Integration-Bulk-CaseIII}}

The second term is a free charge of order $(0,1)$ (see definition \ref{freecharge}). The third term gives us polarization. We rewrite the above as
\bml{\text{term I}&=\sum_{\hat{\ve{x}}'_p\in T(\Box)}l^2\int_{\Box}\int_{-1}^1\left[\dfrac{\tilde{G}(\ve{x},\ve{x}')-\tilde{G}(\ve{x},\hat{\ve{x}}'_p)-\left(l\ve{y}'+hz'\ve{e}_3\right)\cdot \nabla_{\ve{x}'}\tilde{G}(\ve{x},\hat{\ve{x}}'_p)}{|l\ve{y}'+hz'\ve{e}_3|}\right]|\dfrac{l}{h}\ve{y}'+z'\ve{e}_3|\tilde{\rho}_{l,h}(\ve{x}')J(\ve{x}')\dm\ve{y}'\dm z'\\
&+\sum_{\hat{\ve{x}}'_p\in T(\Box)}l^2\tilde{G}(\ve{x},\hat{\ve{x}}'_p)\tilde{q}_{l,h}(\hat{\ve{x}}'_p)J_0(\hat{\ve{x}}'_p)+\sum_{\hat{\ve{x}}'_{p}\in T(\Box)}l^2\nabla_{\ve{x}'}\tilde{G}(\ve{x},\hat{\ve{x}}'_p)\cdot\left(\dfrac{l}{h}\tilde{\ve{p}}_{p,l,h}(\hat{\ve{x}}'_p)+\tilde p_{3,l,h}(\hat{\ve{x}}'_p)\ve{e}_3\right)J_0(\hat{\ve{x}}'_p),\label{Taylor-Integration-Bulk-CaseIII-ii}}
where we defined
\bml{J_0(\hat{\ve{x}}'_p)\tilde{\ve{p}}_{p,l,h}(\hat{\ve{x}}'_p)&=\int_{\Box}\int_{-1}^1 \ve{y}'\tilde{\rho}_{l,h}(\hat{\ve{x}}'_p+l\ve{y}'+hz'\ve{e}_3)J(\hat{\ve{x}}'_p+l\ve{y}'+hz'\ve{e}_3)\dm\ve{y}'\dm z',\\
J_0(\hat{\ve{x}}'_p)\tilde p_{3,l,h}(\hat{\ve{x}}'_p)&=\int_{\Box}\int_{-1}^1z'\tilde{\rho}_{l,h}(\hat{\ve{x}}'_p+l\ve{y}'+hz'\ve{e}_3)J(\hat{\ve{x}}'_p+l\ve{y}'+hz'\ve{e}_3)\dm\ve{y}'\dm z',
}
where $J_0=\textstyle\sqrt{D\psi_0^TD\psi_0}$ is the Jacobian corresponding to the integration on the manifold $\Omega$. Since $\ve{r}\not\in\Omega$, $G(\ve{r},\ve{r}')$ does not contain a singularity. Thus, $G$ has the required smoothness for us to compute derivatives. Noting that ${l}/{h}\to0$ in this limit, the term $({l}/{h})\ve{y}'+z'\ve{e}_3\to z' $ a finite limit. Since $\tilde{G}$ is differentiable, the first term converges to zero from the definition of the derivative. Further, note that $\tilde{\rho}_{l,h}$ and $J$ are bounded and $\textstyle\sum_{T(\Box)} l^2\to\textstyle\int_T\dm\ve{x}'_p$. Thus, the limit of the first term is zero. Taking the limit $h,l\to0$, we get
\beol{\text{term I}\mapsto\int_T\left(\tilde{G}(\ve{x},\ve{x}'_p)\tilde{q}(\ve{x}'_p)+
\dfrac{\partial}{\partial x'_3}\tilde{G}(\ve{x},\ve{x}'_p) \tilde{{p}}_3(\ve{x}'_p)\right)J_0(\hat{\ve{x}}'_p)\dm\ve{x}'_p,\label{BulkLimit-CaseIII}}
where in the limit the tangential polarization $\ve{p}_p$ disappears and $\hat{\ve{x}}_p\mapsto\ve{x}_p$  giving us an integral over $T$. Notice that we cannot integrate out the derivative. The derivative of the Green's function is commonly known as the double-layer potential. This potential arises when we consider the dipole limit for charge surfaces.

Using \eqref{BulkLimit-CaseIII}, we can write the potential as
\bml{\tilde{\Phi}_{0}(\ve{x})&=\int_{T}\left(\tilde{G}(\ve{x},\ve{x}'_p)\tilde{q}(\ve{x}'_p)+\dfrac{\partial}{\partial x'_3}\tilde{G}(\ve{x},\ve{x}'_p)\tilde{{p}}_3(\ve{x}'_p)\right)J_0( {\ve{x}}'_p)\dm\ve{x}'_p}

We can now transform the integral back into real space to get
\beol{\Phi(\ve{r})=\int_{\Omega} {G}(\ve{r},\ve{r}'){q}(\ve{r}')+\dfrac{\partial}{\partial\nu'}G(\ve{r},\ve{r}')p_3(\ve{r}')\dm\ve{r}',\label{Limit-Integral-CaseIII}}
where $\ve{\nu}'$ is the normal to the surface. 

The first part is the free charge, while the second part is the double-layer potential. From the discussion in Appendix \ref{Appendix}, we see that we can rewrite the above as
\bml{-\Delta\Phi_{0}&= q-{p}_3\dfrac{\partial \mathbbm{1}_{\Omega}}{\partial \nu}~~\text{in}~\Scr{R}^3\\
        \Phi_{0}(\ve{r})&=0~~\text{as}~|\ve{r}|\to\infty,\label{PDE-Case3-Homogenized}}
which is \eqref{PDE-Case3-Homogenize}.
\end{proof}

\section{Non-uniqueness of the Polarization Field}\label{Non-Uniqueness}
The polarization, defined as the dipole moment per unit volume, is inherently non-unique. Addressing this non-uniqueness in the polarization distribution requires careful consideration of the arbitrariness in the choice of the unit cell. The formulation of the problem allows for flexibility in selecting the origin \( O \), the lattice \( \scr{L} \), the manifold \( \Omega \) (see discussion on multi-lattice structures in Section \ref{PF}), and the shape of the unit cell \( \Box \). We collectively denote this set of choices by the tuple \( (O, \scr{L}, \Omega, \Box) \).

We consider Theorem \eqref{Theorem-Homo}(2), as it represents the most general case, which, under appropriate conditions, reduces to Theorem \eqref{Theorem-Homo}(1) and (3). Let \(\{\rho_{l,h}\}\) denote a sequence of charge distributions. The free charge, polarization, and total surface charge corresponding to a given choice of parameters \((O, \scr{L}, \Omega, \Box)\) are denoted by \( q:\Omega\to\Scr{R} \), \( \ve{p} : \Omega \to \Scr{R}^3 \), and \( \sigma + \ve{p}_p \cdot \ve{n} : \partial\Omega \to \Scr{R} \), respectively. If, instead, an alternative tuple \( (\utilde{O}, \utilde{\scr{L}}, \utilde{\Omega}, \utilde{\Box}) \) were selected, the corresponding free charge, polarization, and total surface charge would be represented as \( \utilde{q}:\utilde{\Omega}\to\Scr{R} \), \( \utilde{\ve{p}} : \utilde{\Omega} \to \Scr{R}^3 \), and \( \utilde{\sigma} + \utilde{\ve{p}}_p \cdot \utilde{\ve{n}} : \partial\utilde{\Omega} \to \Scr{R} \), respectively. In general, these quantities are expected to differ. 

We assert that, although \( q \), \( \ve{p}_p \), \( p_3 \), and \( \sigma \) are not uniquely defined, the potential obtained from the homogenized problem \eqref{PDE-Case2-Homogenize} remains unique. Here, the polarization has been decomposed into its tangential and normal components. To establish this claim, we observe that the sequence of potentials \( \{\Phi_{l,h}\} \) converges weakly to \( \Phi_0 \). The potentials \( \Phi_{l,h} \) satisfy Poisson’s equation for the charge distribution \( \{\rho_{l,h}\} \). Since the limiting equations are derived directly from the original charge distribution \( \{\rho_{l,h}\} \), the resulting homogenized equation, corresponding to a given choice of domain, lattice, and unit cell, is uniquely determined. Consequently, the associated homogenized potential for this choice is also uniquely defined.

Now, consider an alternative selection of the domain, origin, lattice, and unit cell. The sequence of potentials \(\{\Phi_{l,h}\}\) continues to satisfy Poisson’s equation for the same charge distribution \(\{\rho_{l,h}\}\) while being subject to identical boundary conditions. By the uniqueness of weak convergence, the limiting homogenized potential remains invariant under these choices, thereby establishing that the homogenized potential is independent of the specific selection.

Since \(\Phi_0\) is uniquely determined, it follows from \eqref{PDE-Case2-Homogenize} that the quantities \(\alpha q -\alpha \divS\left(J_0\ve{p}_p\right)/J_0 - \alpha^2p_3{\partial \mathbbm{1}_{\Omega}}/{\partial \nu} \) and \( \alpha\sigma + \alpha\ve{p}_p \cdot \ve{n} \) are uniquely defined and remain independent of the choice of unit cells. Consequently, the non-uniqueness in the tangential component of the polarization is precisely compensated by the surface charge density, ensuring that the resulting potential remains well-defined. For a fixed manifold \( \Omega \), corresponding to a mono-layer 2D material, the normal component of the polarization is uniquely determined.

Most dielectric materials, including piezoelectric and flexoelectric materials, function as insulators and generally do not contain free charge. In such cases, deviations from charge-neutral unit cells are of higher order, ensuring that the free charge of order \((1,0)\) or \((0,1)\) is zero. Consequently, from the preceding analysis, it follows that the expressions $- \alpha\divS\left(J_0\ve{p}_p\right)/J_0 - \alpha^2p_3{\partial \mathbbm{1}_{\Omega}}/{\partial \nu}$ and 
$\alpha(\sigma + \ve{p}_p \cdot \ve{n})$ are uniquely determined and remain independent of the specific choice of unit cells. Furthermore, the influence of curvature is limited to modifying the bound charge term (\( -\divergence~\ve{p} \)). Therefore, the non-uniqueness in the polarization field is effectively compensated by the surface charge distribution, collectively ensuring the uniqueness of the resulting electrostatic potential, associated fields, and energies.

\section*{}

\paragraph*{Acknowledgments.}
    We thank NSF (DMS 2108784, DMS 2342349), AFOSR (MURI FA9550-18-1-0095), and ARO (MURI W911NF-24-2-0184) for financial support.
    Shoham Sen was also partly supported by a Vannevar Bush Faculty Fellowship at the University of Minnesota (PI: R. D. James).
    Kaushik Dayal thanks Air Force Research Laboratory for hosting his visits.

\appendix

\makeatletter
\renewcommand*{\thesection}{\Alph{section}}
\renewcommand*{\thesubsection}{\thesection.\arabic{subsection}}
\renewcommand*{\p@subsection}{}
\renewcommand*{\thesubsubsection}{\thesubsection.\arabic{subsubsection}}
\renewcommand*{\p@subsubsection}{}
\makeatother

\section{The Double-Layer Potential}\label{Appendix}

We investigate the problem of a double-layer potential with the objective of deriving the associated partial differential equation (PDE). We begin by analyzing the case of a flat manifold, which provides a foundational understanding, and then extend our analysis to the more general setting of an arbitrary curved surface.

Consider the following Poisson's equation
\bml{-\Delta\Phi_t(\ve{x})&=\sigma(\ve{x})\mathbbm{1}_{\Omega}(\ve{x}_p)\left(\dfrac{\delta(x_3)-\delta(x_3+t)}{t}\right)~~\text{in}~\Scr{R}^3\\
\Phi_t(\ve{x})&=0~~\text{as}~|\ve{x}|\to\infty,\label{double-layer-Appendix}}
where \(\Omega \subset \Scr{R}^2\), and \(\ve{x} = \ve{x}_p + x_3 \ve{e}_3\), with \(\ve{x}_p \in \Omega\) representing the planar component of \(\ve{x}\) and \(x_3\) representing its normal component. Our focus is on analyzing the behavior of the system as \(t \to 0\).

The physical interpretation of the above result, as discussed in \cite{folland2020introduction}, is that it represents the limit as \(t \to 0\) of the potential induced by a charge distribution with density \(\tfrac{1}{t}\sigma(\ve{x}_p)\) on \(\Omega\), in conjunction with a charge distribution of density \(-\tfrac{1}{t}\sigma(\ve{x}_p)\) on the parallel surface \(\Omega_t = \Omega \oplus \left(-t \ve{e}_3\right)\). This is the dipole limit for a sheet of charge.

The solution for finite $t$ is given by convolution against the Green's function:
\bml{\Phi_t(\ve{x})&=\int_{\Scr{R}^3}G(\ve{x},\ve{x}'){\rho(\ve{x}')}\dm \ve{x}'=\int_{\Omega}G(\ve{x},\ve{x}'_p)\dfrac{\sigma(\ve{x}'_p)}{t}\dm\ve{x}'-\int_{\Omega}G(\ve{x},\ve{x}'_p-t\ve{e}_3)\dfrac{\sigma(\ve{x}'_p)}{t}\dm \ve{x}'\\
&=\int_{\Omega}\left(\dfrac{G(\ve{x},\ve{x}'_p)-G(\ve{x},\ve{x}'_p-t\ve{e}_3)}{t}\right) \sigma(\ve{x}'_p) \dm \ve{x}',}
where for $\ve{x}'\in\Omega_t$, we have set $\ve{x}'=\ve{x}'_p-t\ve{e}_3$. Making this substitution, we get

If we now take the limit as $t\to0$, we get
\beol{\Phi_0(\ve{x})=\int_{\Omega}\lim_{t\to0}\left(\dfrac{G(\ve{x},\ve{x}'_p)-G(\ve{x},\ve{x}'_p-t\ve{e}_3)}{t}\right) \sigma(\ve{x}'_p) \dm \ve{x}'=\int_{\Omega}\dfrac{\partial G(\ve{x},\ve{x}'_p)}{\partial x_3'} \sigma(\ve{x}'_p) \dm \ve{x}'.\label{(90)-Appendix}}

We now take the limit of the right-hand side of \eqref{double-layer-Appendix} to get
\bml{\Delta\Phi_0(\ve{x})&=\sigma(\ve{x}_p)\mathbbm{1}_{\Omega}(\ve{x}_p)\dfrac{\partial \delta(x_3)}{\partial x_3}~~\text{in}~\Scr{R}^3\\
\Phi_0(\ve{x})&=0~~\text{as}~|\ve{x}|\to\infty.\label{double-layer-limit}}

From uniqueness, the solution of \eqref{double-layer-limit} is \eqref{(90)-Appendix}.

We now extend our analysis to the case of general manifolds. Consider the following partial differential equation:
\bml{-\Delta\Phi_t(\ve{r})&=\sigma(\ve{r})\left(\dfrac{\mathbbm{1}_{\Omega}(\ve{r})-\mathbbm{1}_{\Omega_t}(\ve{r})}{t}\right)~~\text{in}~\Scr{R}^3\\
\Phi_t(\ve{r})&=0~~\text{as}~|\ve{r}|\to\infty,\label{double-layer-Appendix-General}}
where $\Omega_t=\{\ve{r}-\ve{\nu}t\vert\ve{r}\in\Omega\}$, where $\ve{\nu}$ is normal to $\Omega$.

The solution to the above is given by 
\beol{\Phi_t(\ve{r})=\dfrac{1}{t}\left(\int_{\Omega}G(\ve{r},\ve{r}'){\sigma(\ve{r}')}\dm\ve{r}'-\int_{\Omega_t}G(\ve{r},\ve{r}'){\sigma(\ve{r}')}\dm\ve{r}'\right)}

Transforming to parameter space using \eqref{MST}, 
\bml{\tilde{\Phi}_t(\ve{x})&=\dfrac{1}{t}\left(\int_{T} \tilde{G}(\ve{x},\ve{x}') {\tilde{\sigma}(\ve{x}')}J(\ve{x}')\dm\ve{x}'-\int_{T_t} \tilde{G}(\ve{x},\ve{x}') {\tilde{\sigma}(\ve{x}')}J(\ve{x}')\dm\ve{x}'\right)\\
&=\int_{T} \left(\dfrac{\tilde{G}(\ve{x},\ve{x}'_p)-\tilde{G}(\ve{x},\ve{x}'_p-t\ve{e}_3)}{t}\right) {\tilde{\sigma}(\ve{x}'_p)}J(\ve{x}'_p)\dm\ve{x}'_p,}
where in the first line, we have set $\tilde{\Phi}_t=(\Phi_t\circ\psi)$, $\tilde{G}=(G\circ\psi)$ and $\tilde{\sigma}=(\sigma\circ\psi)$, while for the second line, we used the substitution $\ve{x}'=\ve{x}'_p+x_3\ve{e}_3$. 

Taking the limit as $t\to0$, we get
\beol{\tilde{\Phi}_0=\int_{T} \dfrac{\partial\tilde{G}(\ve{x},\ve{x}'_p)}{\partial x_3'} {\tilde{\sigma}(\ve{x}'_p)}J(\ve{x}'_p)\dm\ve{x}'_p}

Transforming back into real space, 
\beol{ {\Phi}_0=\int_{\Omega} \dfrac{\partial\tilde{G}(\ve{r},\ve{r}')}{\partial \nu'} {{\sigma}(\ve{r}')}\dm\ve{r}'\label{sol-general-Appendix}}

We now take the limit $t\to0$ of \eqref{double-layer-Appendix-General}
\bml{\Delta\Phi_0(\ve{r})&=\sigma(\ve{r})\dfrac{\partial \mathbbm{1}_{\Omega}}{\partial\nu}~~\text{in}~\Scr{R}^3\\
\Phi_t(\ve{r})&=0~~\text{as}~|\ve{r}|\to\infty,\label{double-layer-Appendix-General-Limit}}

By the uniqueness of solutions, the solution to \eqref{double-layer-Appendix-General-Limit} is given by \eqref{sol-general-Appendix}.

\begin{remark}
    As a consistency check, we observe that setting $\Omega = \omega \times \{0\}$ results in \(\mathbbm{1}_{\Omega} = \mathbbm{1}_{\omega} \mathbbm{1}_{\{0\}} = \mathbbm{1}_{\omega} \delta(x_3)\). This is the same setting as the flat case discussed earlier.
\end{remark}


\newcommand{\etalchar}[1]{$^{#1}$}
\providecommand{\bysame}{\leavevmode\hbox to3em{\hrulefill}\thinspace}
\providecommand{\MR}{\relax\ifhmode\unskip\space\fi MR }
\providecommand{\MRhref}[2]{%
	\href{http://www.ams.org/mathscinet-getitem?mr=#1}{#2}
}
\providecommand{\href}[2]{#2}


\begin{thebibliography}{PdRCL{\etalchar{+}}23}
	
	\bibitem[AA15]{Irene1}
	Amir Abdollahi and Irene Arias, \emph{{Constructive and Destructive Interplay Between Piezoelectricity and Flexoelectricity in Flexural Sensors and Actuators}}, Journal of Applied Mechanics \textbf{82} (2015), no.~12, 121003.
	
	\bibitem[ABC08]{alicandro2008continuum}
	Roberto Alicandro, Andrea Braides, and Marco Cicalese, \emph{Continuum limits of discrete thin films with superlinear growth densities}, Calculus of Variations and Partial Differential Equations \textbf{33} (2008), 267--297.
	
	\bibitem[APM{\etalchar{+}}14]{abdollahi2014computational}
	Amir Abdollahi, Christian Peco, Daniel Millan, Marino Arroyo, and Irene Arias, \emph{Computational evaluation of the flexoelectric effect in dielectric solids}, Journal of Applied Physics \textbf{116} (2014), no.~9, 093502.
	
	\bibitem[AR95]{antman1995nonlinear}
	Stuart~S Antman and Michael Renardy, \emph{Nonlinear problems of elasticity}, SIAM Review \textbf{37} (1995), no.~4, 637.
	
	\bibitem[BBC20]{bach2020discrete}
	Annika Bach, Andrea Braides, and Marco Cicalese, \emph{Discrete-to-continuum limits of multibody systems with bulk and surface long-range interactions}, SIAM Journal on Mathematical Analysis \textbf{52} (2020), no.~4, 3600--3665.
	
	\bibitem[BDM{\etalchar{+}}13]{biswas2013coherent}
	Sushmita Biswas, Jinsong Duan, Krishnamurthy Mahalingam, Dhriti Nepal, Ruth Pachter, Larry Drummy, Dean Brown, and Richard~A Vaia, \emph{Coherent plasmonic engineering in self-assembled reduced symmetry nanostructures}, Plasmonics: Metallic Nanostructures and Their Optical Properties XI, vol. 8809, SPIE, 2013, pp.~180--185.
	
	\bibitem[BNJJ19]{bucsek2019direct}
	Ashley Bucsek, William Nunn, Bharat Jalan, and Richard~D James, \emph{Direct conversion of heat to electricity using first-order phase transformations in ferroelectrics}, Physical Review Applied \textbf{12} (2019), no.~3, 034043.
	
	\bibitem[CAS21]{codony2021transversal}
	David Codony, Irene Arias, and Phanish Suryanarayana, \emph{Transversal flexoelectric coefficient for nanostructures at finite deformations from first principles}, Physical Review Materials \textbf{5} (2021), no.~3, L030801.
	
	\bibitem[CCI{\etalchar{+}}11]{chen2011weak}
	Xian Chen, Shanshan Cao, Teruyuki Ikeda, Vijay Srivastava, G~Jeffrey Snyder, Dominique Schryvers, and Richard~D James, \emph{A weak compatibility condition for precipitation with application to the microstructure of pbte--sb2te3 thermoelectrics}, Acta materialia \textbf{59} (2011), no.~15, 6124--6132.
	
	\bibitem[CDZ{\etalchar{+}}09]{cicalese2009discrete}
	Marco Cicalese, Antonio DeSimone, Caterina~Ida Zeppieri, et~al., \emph{Discrete-to-continuum limits for strain-alignment-coupled systems: Magnetostrictive solids, ferroelectric crystals and nematic elastomers.}, Networks Heterog. Media \textbf{4} (2009), no.~4, 667--708.
	
	\bibitem[CGK24]{cheng2024growth}
	Hsu-Cheng Cheng, Laurent Guin, and Dennis~M Kochmann, \emph{Growth of ferroelectric domain nuclei: Insight from a sharp-interface model}, Journal of the Mechanics and Physics of Solids \textbf{192} (2024), 105810.
	
	\bibitem[CH08]{courant2008methods}
	Richard Courant and David Hilbert, \emph{Methods of mathematical physics, volume 1}, vol.~1, John Wiley \& Sons, 2008.
	
	\bibitem[Che08]{chen2008phase}
	Long-Qing Chen, \emph{Phase-field method of phase transitions/domain structures in ferroelectric thin films: a review}, Journal of the American Ceramic Society \textbf{91} (2008), no.~6, 1835--1844.
	
	\bibitem[FJM02]{friesecke2002theorem}
	Gero Friesecke, Richard~D James, and Stefan M{\"u}ller, \emph{A theorem on geometric rigidity and the derivation of nonlinear plate theory from three-dimensional elasticity}, Communications on Pure and Applied Mathematics: A Journal Issued by the Courant Institute of Mathematical Sciences \textbf{55} (2002), no.~11, 1461--1506.
	
	\bibitem[Fol20]{folland2020introduction}
	Gerald~B Folland, \emph{Introduction to partial differential equations}, Princeton university press, 2020.
	
	\bibitem[GD20]{grasinger2020statistical}
	Matthew Grasinger and Kaushik Dayal, \emph{Statistical mechanical analysis of the electromechanical coupling in an electrically-responsive polymer chain}, Soft Matter \textbf{16} (2020), no.~27, 6265--6284.
	
	\bibitem[GD21]{grasinger2021architected}
	\bysame, \emph{Architected elastomer networks for optimal electromechanical response}, Journal of the Mechanics and Physics of Solids \textbf{146} (2021), 104171.
	
	\bibitem[GDdP22]{grasinger2022statistical}
	Matthew Grasinger, Kaushik Dayal, Gal deBotton, and Prashant~K Purohit, \emph{Statistical mechanics of a dielectric polymer chain in the force ensemble}, Journal of the Mechanics and Physics of Solids \textbf{158} (2022), 104658.
	
	\bibitem[Ger10]{gerstner2010nobel}
	Ed~Gerstner, \emph{Nobel prize 2010: Andre geim \& konstantin novoselov}, Nature Physics \textbf{6} (2010), no.~11, 836--836.
	
	\bibitem[GKK{\etalchar{+}}22]{ghosh2022ferroelectricity}
	Sujoy~Kumar Ghosh, Jinyoung Kim, Minsoo~P Kim, Sangyun Na, Jeonghoon Cho, Jae~Joon Kim, and Hyunhyub Ko, \emph{Ferroelectricity-coupled 2d-mxene-based hierarchically designed high-performance stretchable triboelectric nanogenerator}, ACS nano \textbf{16} (2022), no.~7, 11415--11427.
	
	\bibitem[GMS21]{grasinger2021flexoelectricity}
	Matthew Grasinger, Kosar Mozaffari, and Pradeep Sharma, \emph{Flexoelectricity in soft elastomers and the molecular mechanisms underpinning the design and emergence of giant flexoelectricity}, Proceedings of the National Academy of Sciences \textbf{118} (2021), no.~21, e2102477118.
	
	\bibitem[GRV{\etalchar{+}}20]{glavin2020emerging}
	Nicholas~R Glavin, Rahul Rao, Vikas Varshney, Elisabeth Bianco, Amey Apte, Ajit Roy, Emilie Ringe, and Pulickel~M Ajayan, \emph{Emerging applications of elemental 2d materials}, Advanced Materials \textbf{32} (2020), no.~7, 1904302.
	
	\bibitem[IK21]{indergand2021effect}
	Roman Indergand and Dennis~M Kochmann, \emph{Effect of temperature on domain wall--pore interactions in lead zirconate titanate: A phase-field study}, Applied Physics Letters \textbf{119} (2021), no.~22.
	
	\bibitem[IVNK20]{indergand2020phase}
	Roman Indergand, A~Vidyasagar, Neel Nadkarni, and Dennis~M Kochmann, \emph{A phase-field approach to studying the temperature-dependent ferroelectric response of bulk polycrystalline pzt}, Journal of the Mechanics and Physics of Solids \textbf{144} (2020), 104098.
	
	\bibitem[JBD23]{jha2023discrete}
	Prashant~K Jha, Timothy Breitzman, and Kaushik Dayal, \emph{Discrete-to-continuum limits of long-range electrical interactions in nanostructures}, Archive for Rational Mechanics and Analysis \textbf{247} (2023), no.~2, 29.
	
	\bibitem[JJ21]{james2021conversion}
	Richard~D James and Bharat Jalan, \emph{Conversion of heat to electricity using phase transformations in ferroelectric oxide capacitors}, March~16 2021, US Patent 10,950,777.
	
	\bibitem[JM94]{james1994internal}
	Richard~D James and Stefan M{\"u}ller, \emph{Internal variables and fine-scale oscillations in micromagnetics}, Continuum Mechanics and Thermodynamics \textbf{6} (1994), 291--336.
	
	\bibitem[JMKD23]{jha2023atomic}
	Prashant~K Jha, Jason Marshall, Jaroslaw Knap, and Kaushik Dayal, \emph{Atomic-to-continuum multiscale modeling of defects in crystals with nonlocal electrostatic interactions}, Journal of Applied Mechanics \textbf{90} (2023), no.~2, 021003.
	
	\bibitem[JS10]{jalan2010large}
	Bharat Jalan and Susanne Stemmer, \emph{Large seebeck coefficients and thermoelectric power factor of la-doped srtio3 thin films}, Applied Physics Letters \textbf{97} (2010), no.~4.
	
	\bibitem[JW98]{james1998magnetostriction}
	Richard~D James and Manfred Wuttig, \emph{Magnetostriction of martensite}, Philosophical magazine A \textbf{77} (1998), no.~5, 1273--1299.
	
	\bibitem[KJM{\etalchar{+}}23]{kumbhakar2023prospective}
	Partha Kumbhakar, Jitha~S Jayan, Athira~Sreedevi Madhavikutty, PR~Sreeram, Appukuttan Saritha, Taichi Ito, and Chandra~Sekhar Tiwary, \emph{Prospective applications of two-dimensional materials beyond laboratory frontiers: A review}, IScience \textbf{26} (2023), no.~5.
	
	\bibitem[KS05]{kohn2005another}
	Robert~V Kohn and Valeriy~V Slastikov, \emph{Another thin-film limit of micromagnetics}, Archive for rational mechanics and analysis \textbf{178} (2005), 227--245.
	
	\bibitem[KS16]{krichen2016flexoelectricity}
	Sana Krichen and Pradeep Sharma, \emph{Flexoelectricity: A perspective on an unusual electromechanical coupling}, Journal of Applied Mechanics \textbf{83} (2016), no.~3.
	
	\bibitem[KSV93]{king1993theory}
	RD~King-Smith and David Vanderbilt, \emph{Theory of polarization of crystalline solids}, Physical Review B \textbf{47} (1993), no.~3, 1651.
	
	\bibitem[KTK22]{kannan2022kinetics}
	Vignesh Kannan, Morgan Trassin, and Dennis~M Kochmann, \emph{Kinetics of ferroelectric switching in poled barium titanate ceramics: Effects of electrical cycling rate}, Materialia \textbf{25} (2022), 101553.
	
	\bibitem[LA16]{liu2016localized}
	Zizhuo Liu and Koray Aydin, \emph{Localized surface plasmons in nanostructured monolayer black phosphorus}, Nano letters \textbf{16} (2016), no.~6, 3457--3462.
	
	\bibitem[LDW{\etalchar{+}}24]{liang2024advancements}
	Xu~Liang, Huiting Dong, Yifan Wang, Qianqian Ma, Hongxing Shang, Shuling Hu, and Shengping Shen, \emph{Advancements of flexoelectric materials and their implementations in flexoelectric devices}, Advanced Functional Materials \textbf{34} (2024), no.~51, 2409906.
	
	\bibitem[LPC{\etalchar{+}}17]{lee2017reliable}
	Ju-Hyuck Lee, Jae~Young Park, Eun~Bi Cho, Tae~Yun Kim, Sang~A Han, Tae-Ho Kim, Yanan Liu, Sung~Kyun Kim, Chang~Jae Roh, Hong-Joon Yoon, et~al., \emph{Reliable piezoelectricity in bilayer wse2 for piezoelectric nanogenerators}, Advanced Materials \textbf{29} (2017), no.~29, 1606667.
	
	\bibitem[LS13]{liu2013flexoelectricity}
	LP~Liu and P~Sharma, \emph{Flexoelectricity and thermal fluctuations of lipid bilayer membranes: Renormalization of flexoelectric, dielectric, and elastic properties}, Physical Review E \textbf{87} (2013), no.~3, 032715.
	
	\bibitem[LTY{\etalchar{+}}23]{lin2023recent}
	Yu-Chuan Lin, Riccardo Torsi, Rehan Younas, Christopher~L Hinkle, Albert~F Rigosi, Heather~M Hill, Kunyan Zhang, Shengxi Huang, Christopher~E Shuck, Chen Chen, et~al., \emph{Recent advances in 2d material theory, synthesis, properties, and applications}, ACS nano \textbf{17} (2023), no.~11, 9694--9747.
	
	\bibitem[LZY{\etalchar{+}}16]{li2016origin}
	Fei Li, Shujun Zhang, Tiannan Yang, Zhuo Xu, Nan Zhang, Gang Liu, Jianjun Wang, Jianli Wang, Zhenxiang Cheng, Zuo-Guang Ye, et~al., \emph{The origin of ultrahigh piezoelectricity in relaxor-ferroelectric solid solution crystals}, Nature communications \textbf{7} (2016), no.~1, 13807.
	
	\bibitem[MD14]{marshall2014atomistic}
	Jason Marshall and Kaushik Dayal, \emph{Atomistic-to-continuum multiscale modeling with long-range electrostatic interactions in ionic solids}, Journal of the Mechanics and Physics of Solids \textbf{62} (2014), 137--162.
	
	\bibitem[MK24a]{mathew2024active}
	Anand Mathew and Yashashree Kulkarni, \emph{Active matter as the underpinning agency for extraordinary sensitivity of biological membranes to electric fields}, arXiv preprint arXiv:2412.16319 (2024).
	
	\bibitem[MK24b]{mathew2024electro}
	\bysame, \emph{An electro-chemo-mechanical theory with flexoelectricity: Application to ionic conductivity of soft solid electrolytes}, Journal of Applied Mechanics \textbf{91} (2024), no.~4.
	
	\bibitem[NKPV20]{nepal2020toward}
	Dhriti Nepal, W~Joshua Kennedy, Ruth Pachter, and Richard~A Vaia, \emph{Toward architected nanocomposites: Mxenes and beyond}, ACS nano \textbf{15} (2020), no.~1, 21--28.
	
	\bibitem[NPB{\etalchar{+}}18]{nunn2018frequency}
	William Nunn, Abhinav Prakash, Arghya Bhowmik, Ryan Haislmaier, Jin Yue, Juan~Maria Garcia~Lastra, and Bharat Jalan, \emph{Frequency-and temperature-dependent dielectric response in hybrid molecular beam epitaxy-grown basno3 films}, APL Materials \textbf{6} (2018), no.~6.
	
	\bibitem[NPV12]{nepal2012high}
	Dhriti Nepal, Kyoungweon Park, and Richard~A Vaia, \emph{High-yield assembly of soluble and stable gold nanorod pairs for high-temperature plasmonics}, Small \textbf{8} (2012), no.~7, 1013--1020.
	
	\bibitem[PdRCL{\etalchar{+}}23]{pace20232d}
	Giuseppina Pace, Antonio~Esau del Rio~Castillo, Alessio Lamperti, Simone Lauciello, and Francesco Bonaccorso, \emph{2d materials-based electrochemical triboelectric nanogenerators}, Advanced Materials \textbf{35} (2023), no.~23, 2211037.
	
	\bibitem[PWWC17]{pohlmann2017thermodynamic}
	Hannah Pohlmann, Jian-Jun Wang, Bo~Wang, and Long-Qing Chen, \emph{A thermodynamic potential and the temperature-composition phase diagram for single-crystalline k1-xnaxnbo3 (0<= x<= 0.5)}, Applied Physics Letters \textbf{110} (2017), no.~10.
	
	\bibitem[PWZ{\etalchar{+}}13]{pradel2013piezotronic}
	Ken~C Pradel, Wenzhuo Wu, Yusheng Zhou, Xiaonan Wen, Yong Ding, and Zhong~Lin Wang, \emph{Piezotronic effect in solution-grown p-type zno nanowires and films}, Nano letters \textbf{13} (2013), no.~6, 2647--2653.
	
	\bibitem[RC23]{ross2023thermodynamics}
	Aiden~M Ross and Long-Qing Chen, \emph{Thermodynamics and ferroelectric properties of pb1-xsrxtio3 solid solutions}, Acta Materialia \textbf{261} (2023), 119405.
	
	\bibitem[RV07]{Resta-Vanderbilt}
	Raffaele Resta and David Vanderbilt, \emph{{Theory of Polarization: A Modern Approach}}, {Physics of Ferroelectrics}, {Topics in Applied Physics}, vol. 105, Springer Berlin / Heidelberg, 2007, pp.~31--68.
	
	\bibitem[Sen22]{Sen2022}
	Shoham Sen, \emph{{Nonlocal Dipolar Interactions in Complex Geometries for Quantum Embedding}}, Ph.D. thesis, Carnegie Mellon University, 2022.
	
	\bibitem[Sha94]{sharp1994electrostatic}
	Kim~A Sharp, \emph{Electrostatic interactions in macromolecules}, Current Opinion in Structural Biology \textbf{4} (1994), no.~2, 234--239.
	
	\bibitem[SRD{\etalchar{+}}24]{sen2024multifunctional}
	Dipanjan Sen, Harikrishnan Ravichandran, Mayukh Das, Pranavram Venkatram, Sooho Choo, Shivasheesh Varshney, Zhiyu Zhang, Yongwen Sun, Jay Shah, Shiva Subbulakshmi~Radhakrishnan, et~al., \emph{Multifunctional 2d fets exploiting incipient ferroelectricity in freestanding srtio3 nanomembranes at sub-ambient temperatures}, Nature communications \textbf{15} (2024), no.~1, 10739.
	
	\bibitem[SS23]{surmenev2023influence}
	Roman~A Surmenev and Maria~A Surmeneva, \emph{The influence of the flexoelectric effect on materials properties with the emphasis on photovoltaic and related applications: A review}, Materials Today \textbf{67} (2023), 256--298.
	
	\bibitem[Ste18]{steigmann2018mechanics}
	David~J Steigmann, \emph{Mechanics and physics of lipid bilayers}, The role of mechanics in the study of lipid bilayers (2018), 1--61.
	
	\bibitem[SWBD]{sen2024rigorous}
	Shoham Sen, Yang Wang, Timothy Breitzman, and Kaushik Dayal, \emph{A real-space two-scale analysis of the polarization density in ionic solids}.
	
	\bibitem[SWBD24]{sen2024nonuniqueness}
	\bysame, \emph{Nonuniqueness in defining the polarization: Nonlocal surface charges and the electrostatic, energetic, and transport perspectives}, Journal of the Mechanics and Physics of Solids \textbf{191} (2024), 105743.
	
	\bibitem[Tag86]{tagantsev1986piezoelectricity}
	Alexander~K Tagantsev, \emph{Piezoelectricity and flexoelectricity in crystalline dielectrics}, Physical Review B \textbf{34} (1986), no.~8, 5883.
	
	\bibitem[Tag87]{tagantsev1987pyroelectric}
	Alexandr~Kirillovich Tagantsev, \emph{Pyroelectric, piezoelectric, flexoelectric, and thermal polarization effects in ionic crystals}, Soviet Physics Uspekhi \textbf{30} (1987), no.~7, 588.
	
	\bibitem[TDK{\etalchar{+}}20]{tian2020hexagonal}
	Xiao-Qing Tian, Jing-Yi Duan, Maryam Kiani, Ya-Dong Wei, Naixing Feng, Zhi-Rui Gong, Xiang-Rong Wang, Yu~Du, and Boris~I Yakobson, \emph{Hexagonal layered group iv--vi semiconductors and derivatives: fresh blood of the 2d family}, Nanoscale \textbf{12} (2020), no.~25, 13450--13459.
	
	\bibitem[Tou56]{toupin1956elastic}
	Richard~A Toupin, \emph{The elastic dielectric}, Journal of Rational Mechanics and Analysis \textbf{5} (1956), no.~6, 849--915.
	
	\bibitem[TVV{\etalchar{+}}13]{tagantsev2013origin}
	Alexander~K Tagantsev, K~Vaideeswaran, Sergey~B Vakhrushev, AV~Filimonov, RG~Burkovsky, A~Shaganov, D~Andronikova, AI~Rudskoy, AQR Baron, H~Uchiyama, et~al., \emph{The origin of antiferroelectricity in pbzro3}, Nature communications \textbf{4} (2013), no.~1, 1--8.
	
	\bibitem[TWSK02]{tadmor-EffHamil}
	EB~Tadmor, UV~Waghmare, GS~Smith, and E~Kaxiras, \emph{{Polarization switching in PbTiO3: an ab initio finite element simulation}}, Acta Materialia \textbf{50} (2002), no.~11, 2989--3002.
	
	\bibitem[Van18]{vanderbilt_2018}
	David Vanderbilt, \emph{Berry phases in electronic structure theory: Electric polarization, orbital magnetization and topological insulators}, Cambridge University Press, 2018.
	
	\bibitem[VKS93]{Vanderbilt_surface}
	David Vanderbilt and R.~D. King-Smith, \emph{{Electric polarization as a bulk quantity and its relation to surface charge}}, Phys. Rev. B \textbf{48} (1993), 4442--4455.
	
	\bibitem[WC24]{wang2024theory}
	Bo~Wang and Long-Qing Chen, \emph{Theory of strain phase equilibria and diagrams}, Acta Materialia \textbf{274} (2024), 120025.
	
	\bibitem[WHZ{\etalchar{+}}16]{wang2016subatomic}
	Xuewen Wang, Xuexia He, Hongfei Zhu, Linfeng Sun, Wei Fu, Xingli Wang, Lai~Chee Hoong, Hong Wang, Qingsheng Zeng, Wu~Zhao, et~al., \emph{Subatomic deformation driven by vertical piezoelectricity from cds ultrathin films}, Science Advances \textbf{2} (2016), no.~7, e1600209.
	
	\bibitem[WK17]{wojnar2017linking}
	Charles~S Wojnar and Dennis~M Kochmann, \emph{Linking internal dissipation mechanisms to the effective complex viscoelastic moduli of ferroelectrics}, Journal of Applied Mechanics \textbf{84} (2017), no.~2, 021006.
	
	\bibitem[WLC21]{wang2021inverse}
	Bo~Wang, Fei Li, and Long-Qing Chen, \emph{Inverse domain-size dependence of piezoelectricity in ferroelectric crystals}, Advanced Materials \textbf{33} (2021), no.~51, 2105071.
	
	\bibitem[WWL{\etalchar{+}}14]{wu2014piezoelectricity}
	Wenzhuo Wu, Lei Wang, Yilei Li, Fan Zhang, Long Lin, Simiao Niu, Daniel Chenet, Xian Zhang, Yufeng Hao, Tony~F Heinz, et~al., \emph{Piezoelectricity of single-atomic-layer mos2 for energy conversion and piezotronics}, Nature \textbf{514} (2014), no.~7523, 470--474.
	
	\bibitem[WWW13]{wu2013taxel}
	Wenzhuo Wu, Xiaonan Wen, and Zhong~Lin Wang, \emph{Taxel-addressable matrix of vertical-nanowire piezotronic transistors for active and adaptive tactile imaging}, Science \textbf{340} (2013), no.~6135, 952--957.
	
	\bibitem[WWZ{\etalchar{+}}23]{wang20232d}
	Hao Wang, Yao Wen, Hui Zeng, Ziren Xiong, Yangyuan Tu, Hao Zhu, Ruiqing Cheng, Lei Yin, Jian Jiang, Baoxing Zhai, et~al., \emph{2d ferroic materials for nonvolatile memory applications}, Advanced Materials (2023), 2305044.
	
	\bibitem[WYS19]{wang2019flexoelectricity}
	Binglei Wang, Shengyou Yang, and Pradeep Sharma, \emph{Flexoelectricity as a universal mechanism for energy harvesting from crumpling of thin sheets}, Physical Review B \textbf{100} (2019), no.~3, 035438.
	
	\bibitem[Xia05]{xiao2005influence}
	Yu~Xiao, \emph{The influence of oxygen vacancies on domain patterns in ferroelectric perovskites}, Ph.D. thesis, California Institute of Technology, 2005.
	
	\bibitem[XZF{\etalchar{+}}15]{xue2015influence}
	Fei Xue, Limin Zhang, Xiaolong Feng, Guofeng Hu, Feng~Ru Fan, Xiaonan Wen, Li~Zheng, and Zhong~Lin Wang, \emph{Influence of external electric field on piezotronic effect in zno nanowires}, Nano Research \textbf{8} (2015), 2390--2399.
	
	\bibitem[ZGDY19]{zhang2019room}
	Jun-Jie Zhang, Jie Guan, Shuai Dong, and Boris~I Yakobson, \emph{Room-temperature ferroelectricity in group-iv metal chalcogenide nanowires}, Journal of the American Chemical Society \textbf{141} (2019), no.~38, 15040--15045.
	
	\bibitem[ZLW{\etalchar{+}}15]{zi2015triboelectric}
	Yunlong Zi, Long Lin, Jie Wang, Sihong Wang, Jun Chen, Xing Fan, Po-Kang Yang, Fang Yi, and Zhong~Lin Wang, \emph{Triboelectric--pyroelectric--piezoelectric hybrid cell for high-efficiency energy-harvesting and self-powered sensing}, Advanced materials \textbf{27} (2015), no.~14, 2340--2347.
	
	\bibitem[ZPT{\etalchar{+}}]{zhoudomain}
	Meng-Jun Zhou, Kun Peng, Yueze Tan, Tiannan Yang, Long-Qing Chen, and Ce-Wen Nan, \emph{Domain evolution, dielectric and piezoelectric response in bent freestanding pbtio3 ferroelectric nanowires}, Dielectric and Piezoelectric Response in Bent Freestanding PbTiO3 Ferroelectric Nanowires.
	
	\bibitem[ZWC22]{zorn2022q}
	Jacob~A Zorn, Bo~Wang, and Long-Qing Chen, \emph{Q-pop-thermo: A general-purpose thermodynamics solver for ferroelectric materials}, Computer Physics Communications \textbf{275} (2022), 108302.
	
	\bibitem[ZZCP21]{zhang2021piezotronics}
	Qin Zhang, Shanling Zuo, Ping Chen, and Caofeng Pan, \emph{Piezotronics in two-dimensional materials}, InfoMat \textbf{3} (2021), no.~9, 987--1007.
	
\end{thebibliography}
\end{document}